# An Introduction to Multi-Point Flux (MPFA) and Stress (MPSA) Finite Volume Methods for Thermo-Poroelasticity

*On occasion of the 70th anniversary of Ivar Aavatsmark*

Jan Martin Nordbotten
Eirik Keilegavlen
Department of Mathematics
University of Bergen

## Abstract

We give a unified introduction to the MPFA- and MPSA-type finite volume methods for Darcy flow and poro-elasticity, applicable to general polyhedral grids. This leads to a more systematic perspective of these methods than has been exposed in previous texts, and we therefore refer to this discretization family as the MPxA methods. We apply this MPxA framework to also define a consistent finite-volume discretization of thermo-poro-elasticity. In order to make the exposition accessible to a wide audience, we avoid much of the technical notation which is used in the research literature, and compensate for this by an expanded summary and literature review of the main properties of the MPxA methods. We close the chapter by a section containing applications to problems with complex geometries and non-linear physics.

## 1. Introduction and historical context

The first so-called Multi-Point Finite Volumes methods were developed in the first half of the 1990's, within the context of numerical discretizations for multi-phase flow in geological porous media [1, 2, 3, 4]. In particular, these methods are constructed to solve the so-called pressure equation, which is a second-order elliptic partial differential equation where it is understood that the material parameter may have very low regularity in space. This equation is best presented as a system of first order equations, consisting of balancing the flux $q$ with a source $r$

$$\nabla \cdot q = r \tag{1.1}$$

Complemented by the constitutive law that the flux is derived from a fluid potential $p$:

$$q + \kappa \nabla p = g \tag{1.2}$$

Here $\kappa$ is the material tensor, which is essentially the permeability to flow. The permeability may be both anisotropic and vary strongly as a function of space due to the complex nature of natural rocks. In order to keep the presentation simple, we have simplified terms. Thus it is understood that in applications, the conservation statement is for the mass flux, while the right-hand side of equation (1.2) is the product of the permeability and gravity, etc. For a detailed physical exposition of equations (1.1-1.2), see e.g. [5, 6, 7, 8].

For problems on the form of equations (1.1-1.2), favorable attributes of a numerical discretization method can be summarized as follows (acknowledging that no list of this form is complete):

A) *Flux balance*: An exact local representation of fluid flux balance is considered essential for stability of multi-phase flow simulations. This is made precise in the sense that Stokes' theorem must hold exactly for some volumes $\omega \in \mathcal{T}$, comprising a reasonably fine partitioning $\mathcal{V}$ of the domain:
$$\int_{\partial \omega} q \cdot n \, dS = \int_{\omega} \psi \, dV \tag{1.3}$$

B) *Accuracy on coarse grids*: In geological porous media, the regularity of coefficients is very low, and thus accuracy is to a large extent equated with accurate handling of material discontinuities, in particular when the discontinuities coincide with the boundaries $\partial \omega$.

C) *Flexible grids*: While many early simulation studies were conducted on regular grids, both anisotropic coefficients, as well as complex geological features, motivates discretizations suitable for complex grids.

D) *Symmetric and positive definite discretization matrix*: A symmetric and positive definite (SPD) matrix allows for application of Conjugate Gradient solvers, which have good performance, in particular with respect to memory usage.

E) *Local flux stencils*: The size of the discretization stencil directly impacts both memory usage, but also floating point operations associated with matrix-vector multiplication. A local expression for the flux (as opposed to a post-processed flux), allows the use of automatic differentiation software for constructing the Jacobian for non-linear problems.

F) *Monotonicity of solution*: The continuous problem has the property that for a positive source term $\psi$, and zero-pressure boundary conditions, the pressure $p$ that solves equations (1.1-1.2) is guaranteed to be positive everywhere in the interior of the domain. Monotonicity is closely related to spurious oscillations, which is a major problem for multi-phase simulations.

G) *Accuracy on fine grids*: As the discretization grid is refined, the truncation error of the discrete approximation should vanish, and the discrete approximation should converge to the continuous solution.

It is perhaps intuitive that all these properties cannot be achieved optimally by any linear discretization. By the late 1980s, it was well understood that none of the existing methods at the time achieved all the favorable properties [9]. These were standard Galerkin finite elements (P1-P1 finite elements or similar), Petrov-Galerkin finite elements (P1-P0 finite elements on staggered grids, also known as Control Volume Finite elements), Mixed Finite Elements (lowest-order Raviart-Thomas for flux and P0 for pressure), or Two-Point Finite Volume methods (still the industry standard for practical simulation). We will make a quick summary of the weakness of each of these discretization methods, to better understand the relative advantages (and disadvantages) of the Multi-Point Finite Volume methods.

Galerkin finite elements is perhaps the most common discretization method available in the field of computational mathematics (for an introduction, see text-books [10, 11]). This discretization method is well-suited for simplicial and Cartesian grids, but for more complex grids the definition of the elements becomes more complicated. While Galerkin finite elements have both local stencils as well as lead to SPD matrices, they need post-processing to obtain a local flux balance [12], and are not particularly well suited to discontinuous permeability coefficients [13].

Petrov-Galerkin finite elements, or Control-Volume Finite Elements (CVFE) as we will refer to the method, attempts to improve over the standard finite element methods by introducing a dual grid around each vertex of the primal grid [14]. On this dual grid, piecewise constant test functions are chosen, so that the local Stokes' equation holds exactly. Nevertheless, the pressure solution $p$ is still represented by finite element functions, so the primal grid must still be relatively simple, and no accuracy is gained over finite element methods with respect to discontinuous permeability coefficients. Furthermore, due to the different choice of elements for the trial and solution spaces, the symmetry of the discretization matrix is lost.

Mixed finite elements (MFE) is another way of generalizing finite element methods [15]. In this approach, the first-order structure indicated in equations (1.1-1.2) is retained explicitly, where the pressure is represented as piecewise constant, while the flux is in a relatively simple space whose divergence is piecewise constant (for relatively simple grids this is the lowest-order Raviart-Thomas space, but defining this space becomes non-trivial even for perturbations of Cartesian grids [16, 17]). The mixed-finite element method is accurate for material contrasts, and has an explicit flux balance. On the other hand, it does not immediately lead to an SPD matrix (without hybridization) and has relatively poor monotonicity properties [18].

Two-Point Finite Volume (TPFV) methods in are a sense the simplest methods satisfying the flux balance. The methods consist of imposing equation (1.3) on any polyhedral partition $\mathcal{T}$ of the domain, and then constructing an approximation to $q \cdot n$ using the pressure values in the two neighbors of any face of the polyhedral partitioning. This simplicity leads to a method satisfying all desired properties A)-F) above, and one could ask if it is the perfect method. Unfortunately, the method is indeed too simple – and in contrast to the three preceding methods discussed – the truncation error only vanishes on a quite restrictive class of grids, and thus in general one can observe convergence to the wrong solution (see e.g. [19]).

The above summary gives some impression of the state-of the art when the multi-point methods were developed. As the name suggests, this family of methods attempts to develop a discretization with favorable properties, not by improving on finite element methods, but rather with basis in the TPFV method. More recently, it has been shown how this development ties back to developments also in the finite element literature, a topic that we will return to at several points in later sections of the chapter.

The first multi-point methods were introduced in two independent papers at the ECMOR conference at Røros, Norway in 1994 [1, 2], and an excellent introduction to the Multi-Point Flux Approximation (MPFA), and references to the early literature, can be found by Aavatsmark [20]. However, since that introductory text was written, these methods have seen significant development, both in terms of applicability to complex problems, but also in terms of a maturing of our understanding of the multi-point methods as a general discretization approach. Our goal with this chapter is therefore to provide a contemporary account of these methods. With concrete reference to Aavatsmark [20], the current text covers a consistent treatment of right-hand-side terms in the constitutive laws, more general continuity conditions, discretization of elasticity and poro-elasticity, and a review of the mathematical analysis of these methods. Moreover, our presentation of the method is based on a more abstract construction than in the introduction by Aavatsmark, more suited to general polyhedral grids.

We preempt some of the later discussion by already announcing some of the main features of the multi-point methods. They are developed to have local flux balance, (relatively) small stencils, and be accurate for challenging grids, including polyhedral grids, and handle accurately heterogeneous permeability fields. It has also been shown that the convergence properties of the methods are good, both for smooth and non-smooth data. The cost of these advantages is that the discretization matrix is only symmetric for simplicial grids, although it is in general positive definite. Monotonicity of the discretization holds subject to conditions which are not prohibitively harsh, but still strict enough to affect some realistic cases.

The chapter is subdivided as follows. In the section 2, we will develop the general principles of multi-point finite volume methods, which we refer to as MPxA methods. We will see that these general principles imply a family of methods for elliptic problems with conservation structure. Building on this, we will in section 3 apply the general principles to three concrete problems: First, fluid flow in porous media, as is the classical motivation for these methods, and leads to the MPFA methods. Secondly, to momentum balance in elastic solids, which leads to the so-called Multi-Point Stress Approximation (MPSA) methods. The MPSA methods are naturally suited combined problem of fluid flow in elastically deformable materials, also known as poroelasticity. Moreover, we also consider the case of thermo-poroelasticity, which includes an advective term in addition to the coupling between heat, flow, and deformation. Having developed the discretization methods for these concrete applications, we will review the mathematical and numerical properties of these methods, as has been reported in literature in Section 4, together with applications to real-world data-sets in Section 5.

## 2. Multi-point finite volume methods

We will structure our presentation of the general construction of multi-point finite volume methods in two parts. In section 2.1, we will present the primal grid and the conservation structure, which is common to all finite volume methods. In section 2.2 we will detail the particular choices which give rise to the so-called multi-point finite volume methods. Our goal throughout the exposition is to be both general yet pedagogical. As a result, the presentation, in particular in section 2.2, deviates significantly from the presentation of these methods found in research articles. In section 2.3, we will discuss aspects related to efficient and stable implementation of these methods.

All the derivations in this section are agnostic to the conservation law of interest (mass or momentum). In order to emphasize this generality, we will in this section refer to MPFA or MPSA methods by the generic acronym MPxA. On the other hand, in order to allow for a streamlined presentation, we will present a rather general concept of the so-called O-methods, thereby excluding the less common variants of the MPFA methods, namely the so-called L-, U-, and Z-methods [3, 21, 22, 23].

## 2.1 Finite volume methods

This section gives the basic notion of a finite volume method for a conservation law, following e.g. [24, 25]. As alluded to in the introduction, a conservation law is a statement of the form

$$\frac{d}{dt}\int_\omega u\, dV + \int_{\partial\omega} n \cdot \tau\, dS = \int_\omega r\, dV \qquad (2.1)$$

Given that we have a domain of interest $\Omega \subset \mathbb{R}^n$, the conservation law is interpreted as follows. We are concerned with a conserved quantity $u$ (e.g. mass, momentum or energy) within *any* measurable subdomain $\omega \subset \Omega$ with external normal vector $n$. The conservation law asserts that the accumulation of $u$ witin $\omega$, is determined by a flux field $\tau$, which may represent mass flux, energy flux or stress, and volumetric sources $r$.

As we are concerned with spatial discretization, we will in the remainder of this section disregard the temporal term, and only consider the steady state of equation (2.1). We note that when the variables are sufficiently regular, equation (1.2) and (2.1) are equivalent due to Stokes' theorem. In absence of such regularity, equation (2.1) is a more general statement than (1.2), and this is the motivation for discretizing equation (2.1) directly. Discretization methods that are developed from this viewpoint are known as *finite volume methods*.

In order to construct a numerical method from equation (2.1), we consider a non-overlapping partitioning of $\Omega$ into a finite set of $N$ subdomains $\omega_k \in \mathcal{T}$, for $k = 1 \ldots N$. An example of such a partitioning for $N = 7$ is given by the solid lines in Figure 2.1.

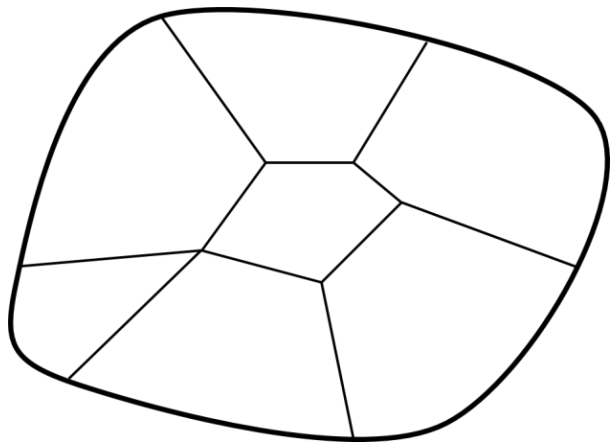

**Figure 2.1:** The domain $\Omega$ shown in thick solid black line together with the finite volume grid $\omega_i$ (thinner solid black lines corresponding to faces between cells). Note that the cells may be polyhedral, and that more than three cells may meet at a vertex.

The subdomains $\omega_k$ are referred to as control volumes, or simpler, *cells*. For any two cells $\omega_{k_1}$ and $\omega_{k_2}$ that are neighbors, in the intersection of their boundaries is measurable, meas $(\partial\omega_{k_1} \cap \partial\omega_{k_2}) > 0$, we refer to this intersection as a *face*, and the collection of faces if denoted $\mathcal{F}$. We extend the definition of a face to also account for intersections with the boundary, such that if meas $(\partial\omega_{k_1} \cap \partial\Omega) > 0$, this also defines a face, and is included in $\mathcal{F}$. In particular, we recognize that all faces of, say, $\omega_k$ is a subset of $\mathcal{F}$, and we denote this subset as $\mathcal{F}_k$. These definitions allow us to rewrite (the steady state of) equation (2.1) as

$$\sum_{\sigma \in \mathcal{F}_k} \int_\sigma n_{\sigma,k} \cdot \tau \, dS = \int_{\omega_k} r \, dV \tag{2.2}$$

Equation (2.2) must hold for any $k$, due to equation (2.1). Moreover, we recognize that it is tempting to define the *normal flux out of $\omega_k$ through $\sigma$* as

$$q_{\sigma,k} \equiv \int_\sigma n_{\sigma,k} \cdot \tau \, dS \tag{2.3}$$

A *finite volume method* is then any method that can be written on the form

$$\sum_{\sigma \in \mathcal{F}_k} q_{\sigma,k} = \int_{\omega_k} r \, dV \qquad \text{for all} \quad \omega_k \in \mathcal{T} \tag{2.4}$$

The finite volume method has local flux balance if for any $\sigma = \partial\omega_{k_1} \cap \partial\omega_{k_2}$, it holds that

$$q_{\sigma,k_1} = -q_{\sigma,k_2} \tag{2.5}$$

Since we will only consider methods with local flux balance, we therefore identify the face flux as the flux from the cell with the lower index, i.e. for $k_1 < k_2$, then we define

$$n_\sigma \equiv n_{\sigma,k_1} \qquad \text{and} \qquad q_\sigma \equiv q_{\sigma,k_1}$$

## 2.2 MPxA Finite volume methods

The basic construction of a finite volume method is agnostic to how the numerical flux field $q_\sigma$ is obtained, and indeed is common for hyperbolic, parabolic and elliptic conservation laws. As stated in the introduction, this chapter deals with methods for problems where there is a proportionality between $q$ and $\nabla u$, as indicated in equation (1.2). In the absence of the time-derivative, such conservation laws are referred to as elliptic, and include Fourier, Fick, Darcy, Hooke and other constitutive laws.

To be precise, we will thus consider constitutive laws on the form

$$\tau = \mathbb{C} \, \nabla u + g \tag{2.6}$$

Here $\mathbb{C}$ is understood to be a local linear operator from the space of functions spanned $\nabla u$, to the space associated with the flux $\tau$. The residual $g$ is in practice derived from a known external force, we will consider it as such. The precise definition of the function spaces depends on the regularity imposed on $u$ and $q$, but also whether one considers scalar of vector equations. As this precision will not be important for introducing the numerical methods, we will omit these details here (for a detailed exposition of the function spaces, see e.g. [26, 27, 24]).

### 2.2.1 Grid structure

The MPxA methods are a family of methods for approximating the normal flux $q_\sigma$ from equation (2.6), based on a core set of foundational principles. To construct an MPxA approximation, additional structure must be introduced relative to the bare-bones finite volume structure given in Section 2.1. In particular, we associate with each cell $\omega_k$ a point $x_k$, which we will refer to as its center. The point $x_k$ should be chosen such that $\omega_k$ is *star-shaped* relative to $x_k$ (this is always possible for simplexes, but may not be possible for non-convex polyhedral). Moreover, we identify that the partition $\mathcal{T}$ gives rise to vertexes of the grid (intersection points of multiple cells), which we will refer to as $\mathcal{V}$. We will denote the subset of $\mathcal{V}$ that are logical vertexes of $\omega_k$ as $\mathcal{V}_k$, such that every point $s \in \mathcal{V}_k$ satisfies $s \in \partial \omega_k$. Conversely, we will denote the subset of $\mathcal{T}$ that meet at a vertrex $s \in \mathcal{V}$ as $\mathcal{T}_s$, such that for every subdomain $\omega_k \in \mathcal{T}_s$, it again holds that $s \in \partial \omega_k$. The definitions of $\mathcal{F}_s$ and $\mathcal{V}_\sigma$, for all $\sigma \in \mathcal{F}$, are analogous.

With the preceding definitions, we introduce a refinement of the finite volume grid structure, as shown in Figure 2.2. First, we refine the faces of the grid as follows: Let every face $\sigma \in \mathcal{F}$ be partitioned into *subfaces* $\tilde{\sigma} \in \mathcal{S}_\sigma$, such that each subface contains exactly one vertex of $\sigma$. Thus, if the set of all subfaces is denoted $\mathcal{S}$, then for any pair of a face $\sigma \in \mathcal{F}$ and a vertex $s \in \mathcal{V}_\sigma$, there is a unique element of $\mathcal{S}_{\sigma,s}$. Extending the notational convention above, we denote the subfaces of $\omega_k$ meeting at a vertex $s \in \mathcal{V}_k$ as $\mathcal{S}_{k,s}$.

We introduce the following definition of a dual grid: For each vertex $s \in \mathcal{V}$, let the dual cell $\omega_s^* \in \mathcal{T}_s^*$ be defined such that subfaces in $\tilde{\sigma} \in \mathcal{S}_s$ are contained in $\omega_s^*$, and the cell-centers $x_k$ of the cells $\omega_k \in \mathcal{T}_s$ are on the boundary of the dual cell, i.e. $x_k \in \partial \omega_s^*$. Finally, let the dual cells be a non-overlapping partitioning of the domain $\Omega$. The intersection of the primal and dual grids creates an even finer grid $\tilde{\mathcal{T}}$, elements of which are uniquely defined by a cell and a corner. Thus, the subcell $\tilde{\omega}_{k,s} \in \tilde{\mathcal{T}}$ is defined as $\tilde{\omega}_{k,s} = \omega_k \cap \omega_s^*$. Again, we retain the same conventions on subscripts, so that in particular $\tilde{\mathcal{T}}_s$ are the subcells adjacent to the corner $s$.

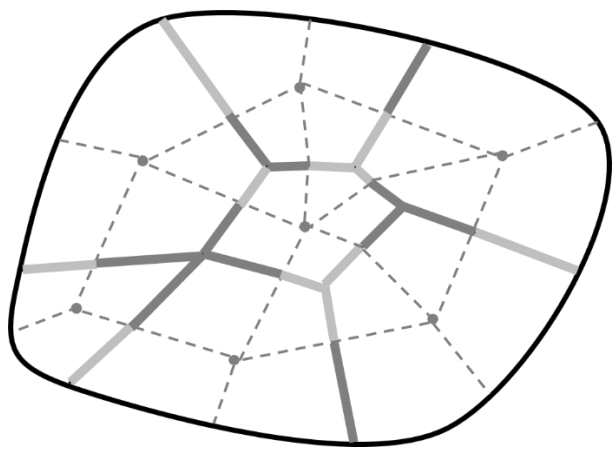

**Figure 2.2:** This figure provides an illustration of the extra grid structure used for MPxA methods relative to the basic finite volume grid shown in Figure 2.1. The division of faces into subfaces is indicated by two nuances of grey, while cell centers are indicated by dots. The dual grid is indicated by dashed lines, and the subgrid is thus obtained as the quadrilaterals having two dashed and two solid-grey boundaries.

### 2.2.2 Approximation spaces

The MPxA approximations do not attempt to construct the numerical flux $q_\sigma$ over a face $\sigma \in \mathcal{F}$ directly, but instead construct approximations over the subfaces $\tilde{\sigma}$. The subface normal fluxes $\tilde{q}_{\tilde{\sigma}}$ are then subsequently assembled such that

$$q_\sigma = \sum_{\tilde{\sigma} \in \mathcal{S}} \tilde{q}_{\tilde{\sigma}} \tag{2.7}$$

The MPxA methods for the fluxes over the subfaces $\tilde{\sigma} \in \mathcal{S}_s$ are based on the following common seven ingredients:

i) A linear approximation $u_{k,s}(x)$ to the potential field within each subcell $\widetilde{\omega}_{k,s} \in \widetilde{\mathcal{T}}$.
ii) A constant approximation $\tau_{k,s}$ to the flux field within each subcell $\widetilde{\omega}_{k,s} \in \widetilde{\mathcal{T}}$.
iii) A constant approximation $g_k$ to the external force field within each cell $\omega_k \in \mathcal{T}$.
iv) A relation between the potential fields $u_{k,s}(x)$, flux fields $\tau_{k,s}$, and force field $g_k$, consistent with (2.6).
v) Local flux balance over each subface $\tilde{\sigma} \in \mathcal{S}_s$ in the sense of (2.3) and (2.5):
$$\int_{\tilde{\sigma}} \tau_{k_1,s} \cdot n_{\tilde{\sigma}} \, dS = \int_{\tilde{\sigma}} \tau_{k_2,s} \cdot n_{\tilde{\sigma}} \, dS \equiv \tilde{q}_{\tilde{\sigma}} \tag{2.8}$$
vi) Continuity between the linear potential field approximations at the cell centers $x_k$, i.e. for a given $\omega_k \in \mathcal{T}$, and any $s_1, s_2 \in \mathcal{V}_k$
$$u_{k,s_1}(x) = u_{k,s_2}(x) \equiv u_k \tag{2.9}$$
vii) A minimization of a quadratic penalty function $\mathcal{M}_s$, measuring the discontinuity of the linear potential fields across subfaces $\tilde{\sigma} \in \mathcal{S}_s$. The precise choice of the penalty function used in this minimization gives rise to variations in the method, as detailed in the next section.

If one of the subfaces $\tilde{\sigma} \in \mathcal{S}_s$ is on the boundary of the domain $\Omega$, additional conditions apply, as discussed in Section 2.2.4. We emphasize that the method formulation given below applies independent of the boundary condition assigned.

The core MPxA ingredients are perhaps most intuitively understood by the following interpretation: The potential field is piecewise linear function relative to the fine grid obtained from the intersection of the primal and dual grid, with a minimum of continuity imposed in order to allow for a compromise between a consistent discretization and flexible grids, while always allowing for a static condensation in terms of cell-center potentials alone. The numerical flux is a derived quantity from the potential field.

The continuity requirements on the piecewise linear potential field are chosen to have a very particular structure. With reference to figure 2.2, continuity conditions are essentially imposed at cell centers (points in the figure), and across subfaces (various thick grey lines). Thus *no continuity is explicitly enforced over the edges of the dual grid*. This observation justifies the claim that all degrees of freedom in the construction can be locally eliminated with respect to the cell-center potentials $u_k$. This is the key to an efficient numerical implementation, and also implies that the resulting discretization matrix has minimum size (potential variables in the cell centers).

To make the above claims more precise, we now introduce discrete operators that allow for an efficient presentation of MPxA methods in general. A more detailed discussion with focus on implementation in provided in Section 2.3, while application of the general framework, that is, identification of the discrete operators for specific equations is considered in Section 3

We denote vectors of variables by bold letters, and matrixes by capitals. First, let the finite volume method, equation (2.4), be represented in matrix form in terms of a divergence matrix $\boldsymbol{D}$ (simply a summation over fluxes, accounting for sign convention). Similarly, we represent the summation over subfluxes, weighted by the area of the subcells, as defined in equation (2.7) as $\boldsymbol{\Sigma}_{\mathcal{F}}$. Then equation (2.4) and (2.7) are equivalent to

$$\boldsymbol{D}\boldsymbol{q} = \boldsymbol{D}\boldsymbol{\Sigma}_{\mathcal{F}}\widetilde{\boldsymbol{q}} = \boldsymbol{r} \tag{2.10}$$

Furthermore, let the vector of cell center potentials be denoted $\boldsymbol{u}$, and the vector containing the degrees of freedom for the linear pressure variations in each subcell $\widetilde{\boldsymbol{u}}$. We denote the operator that extracts $\boldsymbol{u}$ from $\widetilde{\boldsymbol{u}}$ as $\boldsymbol{E}$, such that

$$\boldsymbol{u} = \boldsymbol{E}\widetilde{\boldsymbol{u}} \tag{2.11}$$

The (continuous) gradient induces a map from $\widetilde{\boldsymbol{u}}$ to piece-wise constant vector fields on each subcell, and we denote the matrix representation of this map as $\boldsymbol{G}$. We denote the discrete constitutive law by the matrix $\boldsymbol{B}$, such that equation (2.6) becomes

$$\boldsymbol{\tau} = \boldsymbol{B}\boldsymbol{G}\widetilde{\boldsymbol{u}} + \boldsymbol{E}^*\boldsymbol{g} \tag{2.12}$$

Here $\boldsymbol{E}^*$ is the matrix that maps cell values to the individual subcells (in a sense dual to $\boldsymbol{E}$).

We furthermore denote by $\boldsymbol{F}$ the matrix that extracts normal subface normal fluxes $\widetilde{\boldsymbol{q}}$ from the fluxes $\boldsymbol{\tau}$ on the side of the face with the *lower* index (similar to definition used in (2.5b)), and conversely we denote by $\widehat{\boldsymbol{F}}$ the matrix that extracts normal subface normal fluxes $\widetilde{\boldsymbol{q}}$ from the fluxes $\boldsymbol{\tau}$ on the side of the face with the *higher* index. The flux balance and the definition of the subface fluxes is summarized in matrix form as

$$\boldsymbol{F}\boldsymbol{\tau} = \widehat{\boldsymbol{F}}\boldsymbol{\tau} \tag{2.13}$$

$$\widetilde{\boldsymbol{q}} = \boldsymbol{F}\boldsymbol{\tau} \tag{2.14}$$

Finally, let the penalty function $\mathcal{M}(\widetilde{\boldsymbol{u}}) = \sum_{s\in\mathcal{V}} \mathcal{M}_s(\widetilde{\boldsymbol{u}})$ be the (still quadratic) measure of the total discontinuity of $\widetilde{\boldsymbol{u}}$ across subfaces $\tilde{\sigma} \in \mathcal{S}$.

Then the MPxA method can then be explicitly defined as follows:

<u>Definition 2.1 (generalized MPxA (global))</u>: Let $\mathcal{N}_{u,g}$ be the null-space of the constraints given in equation (2.11 - 2.13), subject to a given potential $\boldsymbol{u}$ and an external field $\boldsymbol{g}$. Then the MPxA method is defined by the pair $(\widetilde{\boldsymbol{u}}, \boldsymbol{\tau}) \in \mathcal{N}_{u,g}$ such that

$$(\widetilde{\boldsymbol{u}}, \boldsymbol{\tau}) = \arg\min_{(\widetilde{\boldsymbol{u}}', \boldsymbol{\tau}')\in\mathcal{N}_{u,g}} \mathcal{M}(\widetilde{\boldsymbol{u}}') \tag{2.15}$$

and the numerical flux is defined as $\boldsymbol{q} = \boldsymbol{Q}_u\boldsymbol{u} + \boldsymbol{Q}_g\boldsymbol{g} \equiv \boldsymbol{\Sigma}_{\mathcal{F}}\boldsymbol{F}(\boldsymbol{B}\boldsymbol{G}\widetilde{\boldsymbol{u}} + \boldsymbol{E}^*\boldsymbol{g})$.

We emphasize that $\mathcal{M}$ is a sum of local quadratic measures $\mathcal{M}_s$ for each $s \in \mathcal{V}$, and moreover that the constraints (2.11-2.13) are all local expressions relative to subcells $\tilde{\omega} \in \tilde{\mathcal{T}}_s$. That is to say that the matrices $\boldsymbol{B}, \boldsymbol{E}, \boldsymbol{E}^*, \boldsymbol{F}, \widehat{\boldsymbol{F}}$ and $\boldsymbol{G}$ can all be written as a sum of local matrices for each vertex, e.g. $\boldsymbol{B} = \sum_{s \in \mathcal{V}} \boldsymbol{B}_s$, where the local matrices such as $\boldsymbol{B}_s$ are in terms of degrees of freedom only associated with the subcells in $\tilde{\mathcal{T}}_s$. This gives rise to the local formulation of MPxA, which is defined as

<u>Definition 2.2 (generalized MPxA (local)):</u> For a vertex $s \in \mathcal{V}$, and for a given potential $\boldsymbol{u}$, let $\mathcal{N}_{u,g,s}$ be the null-space of the constraints

$$\boldsymbol{u} = \boldsymbol{E}_s \tilde{\boldsymbol{u}}_s, \qquad \boldsymbol{\tau}_s = \boldsymbol{B}_s \boldsymbol{G}_s \tilde{\boldsymbol{u}}_s + \boldsymbol{E}_s^* \boldsymbol{g}, \quad \text{and} \quad \boldsymbol{F}_s \boldsymbol{\tau}_s = \widehat{\boldsymbol{F}}_s \boldsymbol{\tau}_s \qquad (2.16)$$

in terms of the local degrees of freedom $\tilde{\boldsymbol{u}}_s$ and $\boldsymbol{\tau}_s$ on subcells in $\tilde{\mathcal{T}}_s$. Then the local problem for the MPxA method is defined by the pair $(\tilde{\boldsymbol{u}}_s, \boldsymbol{\tau}_s) \in \mathcal{N}_{u,s}$ such that

$$(\tilde{\boldsymbol{u}}_s, \boldsymbol{\tau}_s) = \arg \min_{(\tilde{\boldsymbol{u}}_s', \boldsymbol{\tau}_s') \in \mathcal{N}_{u,g,s}} \mathcal{M}_s(\tilde{\boldsymbol{u}}_s') \qquad (2.17)$$

and the numerical flux is assembled as $\boldsymbol{q} = \boldsymbol{Q}_u \boldsymbol{u} + \boldsymbol{Q}_g \boldsymbol{g} \equiv \Sigma_{\mathcal{F}} \sum_{s \in \mathcal{V}} \boldsymbol{F}_s (\boldsymbol{B}_s \boldsymbol{G}_s \tilde{\boldsymbol{u}}_s + \boldsymbol{E}_s^* \boldsymbol{g})$.

The local formulation of MPxA is clearly equivalent to the global formulation. As a consequence, the minimization problems (2.15) are local linear saddle-point problems of modest size for each vertex of the grid, and can be solved efficiently (and in parallel, if desired), using any standard explicit linear solver. We shall return to the structure of the local problems in Section 2.3.3.

Since the minimization problem is quadratic, the numerical flux is a linear function of the potential $\boldsymbol{u}$. The MPxA finite volume discretization matrix is obtained by combining the numerical flux and the finite volume method, equation (2.10), to obtain the linear system

$$\boxed{\boldsymbol{D}\boldsymbol{Q}_u \boldsymbol{u} = \boldsymbol{r} - \boldsymbol{D}\boldsymbol{Q}_g \boldsymbol{g}} \qquad (2.18)$$

We will return to the properties of the matrix $\boldsymbol{D}\boldsymbol{Q}_u$ in Section 4.

### 2.2.3 Penalty functions

An attractive feature of the MPxA methods is that the penalty functions $\mathcal{M}_s$ used in minimization problems (2.15) can be chosen to enhance various properties of the MPxA methods.

The natural starting point for developing quadratic minimization problems to penalize the discontinuities in the linear pressure approximation, is to consider the norm of the discontinuities across subfaces [28]. Thus, for every subface $\tilde{\sigma} \in \mathcal{S}$, we define the penalty function

$$\mathcal{M}_{\tilde{\sigma}}(\boldsymbol{u}) \equiv \int_{\tilde{\sigma}} \left( u_{k_1,s}(x) - u_{k_2,s}(x) \right)^2 dS \qquad (2.19)$$

as previously, $k_1$, $k_2$ and $s$ are the indexes such that $\tilde{\omega}_{k_1,s}$ and $\tilde{\omega}_{k_1,s}$ are the two subcells sharing the subface $\tilde{\sigma}$. Any positive linear combination of the subface discontinuity measure will be a new measure of discontinuity, and thus it is follows that for any vertex, we make the natural definition

$$\mathcal{M}_s(\boldsymbol{u}) \equiv \sum_{\tilde{\sigma} \in \mathcal{S}_s} c_{\tilde{\sigma}} \mathcal{M}_{\tilde{\sigma}}(\boldsymbol{u}) \qquad (2.20)$$

The weights $c_{\tilde{\sigma}}$ can in principle be chosen to optimize the method, although the simple choice $c_{\tilde{\sigma}} = 1$ appears sufficient in practice.

Since the potentials $u_{k,s}(x)$ are approximated as linear, the integral in equation (2.19) is a quadratic function on the subface $\tilde{\sigma}$, and can be exactly evaluated using only a low number of quadrature points (two in 2D and four in 3D). The majority of MPxA literature simplify the minimization problem further, and consider only a single quadrature point. We will for historic reasons denote this minimization with the Greek letter $\eta$, and introduce the simplified penalty functions

$$\mathcal{M}_{\tilde{\sigma}}^{\eta}(\boldsymbol{u}) \equiv \left(u_{k_1,s}(x_{\tilde{\sigma}}^{\eta}) - u_{k_2,s}(x_{\tilde{\sigma}}^{\eta})\right)^2 \tag{2.21}$$

The definition of the simplified penalty functions is completed by specifying the points $x_{\tilde{\sigma}}^{\eta}$. The common choice is obtained if the face $\sigma$ subdivided into subfaces relative to a central point $x_\sigma$. Then let $\eta \in [0,1)$, and define

$$x_{\tilde{\sigma}}^{\eta} = x_\sigma + \eta \frac{x_s - x_\sigma}{|x_s - x_\sigma|} \tag{2.22}$$

In this expression, we have used $x_s$ to denote the coordinate of vertex $s$. Given $\mathcal{M}_{\tilde{\sigma}}^{\eta}(\boldsymbol{u})$, the full expression for minimization $\mathcal{M}_s^{\eta}(\boldsymbol{u})$ is defined analogously to equation (2.20).

The main advantage of the simplified penalty functions $\mathcal{M}_s^{\eta}$, is that it can be shown that for many common grid types (all grids in 2D, and e.g. Cartesian or simplicial grids in 3D, but not grids containing pyramids), the optimal value of the minimization problems (2.17) is indeed 0. That is to say, that the minimization problem can be omitted, and be replaced by the direct condition that

$$\mathcal{M}_{\tilde{\sigma}}^{\eta}(\boldsymbol{u}) = 0 \tag{2.23}$$

for all subfaces $\tilde{\sigma} \in \mathcal{S}$. When this condition holds, the pressure is indeed continuous across the subface exactly at the point $x_{\tilde{\sigma}}^{\eta}$, and this point is then referred to as a continuity point in the literature [20].

Two particular choices of $x_{\tilde{\sigma}}^{\eta}$ are particularly appealing and common in practice: $\eta = 0$ leads to a simple method that has the best monotonicity properties on quadrilaterals [29]. $\eta = \frac{1}{3}$ leads to a method which has a symmetric discretization method on simplexes [30, 31]. Another possible choice is to take $\eta = \frac{1}{2}$, which gives a high-order method on smooth problems on quadrilaterals [32]. We will return to this topic in more detail in Section 4 of the chapter.

## 2.3 Implementation aspects

To further explore the approximation properties and implementation of the MPxA methods, it is instructive to consider the local problem 2.2 in some more detail. As discussed above, the approximation spaces on the subcells are not rich enough to allow full continuity over subfaces $s \in \mathcal{S}_s$, thus MPxA can be interpreted as a discontinuous Galerkin method with a particular set of continuity constraints on potentials and normal fluxes over the subfaces. Critical for efficiency and implementation, we exploit the two-scale approximation, in that the (fine scale) degrees of freedom associated with the potential gradients on the subcells can be eliminated by static condensation around each vertex $s$. This leaves a method where only the (coarse scale) cell center degrees of freedom enter into the global problem.

### 2.3.1 Local minimization problem

To be concrete, we make the choice of representing the linear potential field, $u_{k,s}$ in a subcell by its value in the cell center, $u_k$, and the (constant) components of its gradient, which we denote $h_{k,s}$. To understand an efficient implementation of the local linear system set in Definition 2.2, it is instructive to discuss the size of the matrices that form the problem. To that end, let $n_f = |\mathcal{S}_s|$ be the number of subfaces meeting in $s$, and similarly $n_c = |\tilde{\mathcal{T}}_s|$ be the number of cells that has $s$ as vertex. The number of faces with Neumann and Dirichlet boundary conditions are denoted $n_{\tilde{\sigma},N}$ and $n_{\tilde{\sigma},D}$, respectively. Let $d$ be the dimension of the potential field $u$, this will be 1 for scalar equations and the spatial dimension $n$ for vector equations, and $n_{\tilde{\sigma},q}$ be the number of quadrature points on subface $\tilde{\sigma}$.

As degrees of freedom in the local linear system, we use the cell center potentials $\boldsymbol{u}_s$ in $\mathcal{T}_s$ and the components of the gradients in the subcells $\tilde{\mathcal{T}}_s$, represented by $\boldsymbol{h}_s$, so that the full vector of local unknowns is $\tilde{\boldsymbol{u}}_s = (\boldsymbol{u}_s, \boldsymbol{h}_s)^T$. This representation has the advantage that the matrices $\boldsymbol{E}$ and $\boldsymbol{G}$ take the particularly simple form

$$\boldsymbol{E}_s \tilde{\boldsymbol{u}}_s = (\boldsymbol{I} \quad \boldsymbol{0}) \tilde{\boldsymbol{u}}_s = \boldsymbol{u}_s$$

$$\boldsymbol{G}_s \tilde{\boldsymbol{u}}_s = (\boldsymbol{0} \quad \boldsymbol{I}) \tilde{\boldsymbol{u}}_s = \boldsymbol{h}_s$$

The flux field can therefore be recovered directly from $\boldsymbol{h}_s$ using (2.12), where we see that the (local) matrix that contains the constitutive law, $\boldsymbol{B}_s$ is of size $(n_c \cdot d \cdot n) \times (n_c \cdot d \cdot n)$.

The computation of subface normal fluxes is split into internal and boundary faces: The internal faces are covered by the matrices $\boldsymbol{F}_s^I$ and $\widehat{\boldsymbol{F}}_s^I$, both of size $\left(d \cdot (n_f - n_{\tilde{\sigma},N} - n_{\tilde{\sigma},D})\right) \times (n_c \cdot d \cdot n)$, which represent the multiplication by subface normal vectors of the flux on the neighboring cells of lower and higher index, respectively. The normal flux over faces with Neumann and Dirichlet boundary conditions is computed from the matrices $\boldsymbol{F}_s^N$ and $\boldsymbol{F}_s^D$, of size $(d \cdot n_{\tilde{\sigma},N}) \times (n_c \cdot d \cdot n)$ and $(d \cdot n_{\tilde{\sigma},D}) \times (n_c \cdot d \cdot n)$, respectively. The evaluation of the potential at the subface quadrature points is similarly split: For internal subfaces, the matrix $\boldsymbol{M}_s^I$ of size $\left(d \cdot \sum_{\tilde{\sigma} \in \mathcal{S}_s^i} n_{\tilde{\sigma},q}\right) \times (n_c \cdot d \cdot n)$ has elements composed of the distance from cell centers to subface quadrature points; here $\mathcal{S}_s^i$ denotes the internal subfaces of vertex $s$. $\widehat{\boldsymbol{M}}_s^I$ is the corresponding matrix for neighboring cells of higher index, while $\boldsymbol{M}_s^D$ and $\boldsymbol{M}_s^N$ are assigned for subfaces with Dirichlet and Neumann boundary conditions, respectively. Finally, we similarly define the matrix $\widehat{\boldsymbol{E}}_s^*$ relative to $\boldsymbol{E}_s^*$, and the internal and boundary components.

With the above definitions, the minimization problem is stated in terms of the local variables as

$$\mathcal{M}_s(\tilde{\boldsymbol{u}}_s) = \mathcal{M}_s(\boldsymbol{u}_s, \boldsymbol{h}_s) = \left\| \Sigma \left( (\boldsymbol{E}_s^{*,I} - \widehat{\boldsymbol{E}}_s^{*,I}) \boldsymbol{u}_s + (\boldsymbol{M}_s^I - \widehat{\boldsymbol{M}}_s^I) \boldsymbol{h}_s \right) \right\|^2 \tag{2.24}$$

with $\Sigma$ a diagonal matrix containing the quadrature weights for the integrals. The minimization problem is subject to the constraints that $\boldsymbol{u}_s$ and $\boldsymbol{g}$ are given, as well as

$$\boldsymbol{\tau}_s = \boldsymbol{B}_s \boldsymbol{h}_s + \boldsymbol{E}_s^* \boldsymbol{g}, \qquad \boldsymbol{F}_s \boldsymbol{\tau}_s = \widehat{\boldsymbol{F}}_s \boldsymbol{\tau}_s, \qquad \boldsymbol{F}_s^N \boldsymbol{\tau}_s = \boldsymbol{q}_s^N, \qquad \boldsymbol{E}_s^{*,D} \boldsymbol{u}_s + \boldsymbol{M}_s^D \boldsymbol{h}_s = \boldsymbol{u}_s^D \tag{2.25}$$

Here, we have retained the explicit dependency of the constraints on $\tau_s$, and introduced the Neumann and Dirichlet conditions as $q_s^N$ and $u_s^D$, respectively. These conditions are void if none of the subfaces around vertex $s$ are located on the domain boundary. We have assumed that the number of quadrature points on subfaces with Dirichlet conditions is sufficiently low for the relevant constraint to be fulfilled exactly; as an alternative, the condition $M_s^D u_s = u_s^D$ can be incorporated into the minimization problem. While the matrices $F_s^D$ and $M_s^N$ are not used in the method constructions, they can be used in post-processing to calculate the normal flux through Dirichlet faces and potential on Neumann faces, respectively.

We pause to consider the size of the local problem (2.24)-(2.25), and specifically compare the number of degrees of freedom and constrains. We limit ourselves here to internal vertexes, similar reasoning applies to vertexes on the domain boundary. The number of subcell gradient degrees of freedom is $n_k \cdot d \cdot n$, while there are $n_{\tilde{\sigma}} \cdot d$ equations for flux continuity, and $\left(d \cdot \sum_{\tilde{\sigma} \in S_s} n_{\tilde{\sigma},q}\right)$ quadrature points for potential continuity. If $n = 2$, $n_{\tilde{\sigma}} = n_c$ independent of the cell type thus if each subfaces is assigned a single quadrature point, $n_{\tilde{\sigma},q} = 1$, then the number of equations equals the number of gradient unknowns. In 3d, the situation is more nuanced: For simplex and logically Cartesian grids, $\frac{n_{\tilde{\sigma}}}{n_c} = \frac{3}{2}$, thus with a single quadrature point on each subface, the number of equations and gradient unknowns still match. For general cell shapes, notably pyramids, this is no longer the case, and the method with a single quadrature point fails.

### 2.3.2 Expression in terms of coarse degrees of freedom

To arrive at a discretization in terms of the coarse scale, cell center, degrees of freedom, the next step is to eliminate the subcell gradients $h_s$. Which approach is practical here depends on the size of the minimization problem. When the number of equations and gradient unknowns match, it turns out that the value of the minimization problem is in fact zero, and the problem can be formulated as a linear system on the form (with $\tau_s$ eliminated)

$$\begin{pmatrix} M_s^I - \widehat{M}_s^I \\ F_s B_s - \widehat{F}_s B_s \\ F_s^N B_s \\ M_s^D \end{pmatrix} h_s = - \begin{pmatrix} E_s^{*,I} - \widehat{E}_s^{*,I} \\ 0 \\ 0 \\ E_s^{*,D} \end{pmatrix} u_s + \begin{pmatrix} 0 \\ (E_s^{*\prime} - \widehat{E}_s^*) g_s \\ q_s^N \\ u_s^D \end{pmatrix} \quad (2.26)$$

This system can be solved to express $h_s$ as a linear function of $u_s$ and the right-hand side terms. We write the respective solutions as $h_s^u = S_u u_s$, $h_s^g = S_g g_s$, $h_s^N = S_N q_s^N$ and $h_s^D = S_D u_s^D$, where the matrices $S_*$ are computed from the left and right hand sides of (2.26). The solvability of (2.26) depends on the grid types and problem under consideration, and problems can arise in special cases such as some non-matching grids in 3D. However, for regular grids (simplicial and Cartesian) for the problems considered in Section 3 (with the exception of subsection 3.2.1), equation (2.26) is solvable.

The more general case with multiple quadrature points leads to a true minimization problem, and thus resolves many of the cases where (2.26) is not suitable. Since this is a quadratic minimization problem with linear constraints, the minimum can be found in the standard way as the solution to a linear system of equation obtained via a Lagrange multiplier vector $\lambda$. For completeness, we state this system for internal cells (i.e. with no boundary cells), for which the constrained system is:

$$\begin{pmatrix} (M_s^I - \widehat{M}_s^I)^T \Sigma^T \Sigma (M_s^I - \widehat{M}_s^I) & (F_s B_s - \widehat{F}_s B_s)^T \\ F_s B_s - \widehat{F}_s B_s & 0 \end{pmatrix} \begin{pmatrix} h_s \\ \lambda \end{pmatrix}$$
$$= -\begin{pmatrix} 2(M_s^I - \widehat{M}_s^I)^T \Sigma^T \Sigma (E_s^{*,I} - \widehat{E}_s^{*,I}) \\ 0 \end{pmatrix} u_s + \begin{pmatrix} 0 \\ (E_s^{*,\prime} - \widehat{E}_s^*) g_s \end{pmatrix}$$

(2.27)

Thus while the local system is still linear, and has the advantage that unique solvability can be proved for many classes of grids due to the relationship to the minimization problem [28]. However, equation (2.27) contains $n_{\tilde{\sigma}} \cdot d$ extra unknowns, corresponding to the Lagrange multipliers for flux continuity.

### 2.3.3 Computational cost

We make the following comments on the local problems: Their size depends on the number of cells that share the vertex $s$, the dimension of the potential field, and the number of quadrature points assigned on the subfaces. In practice, 12 gradients are eliminated for scalar problems on a Cartesian 2d grid with a single quadrature points, up to on the order of a few hundred degrees of freedom for vector problems on 3d unstructured grids with multiple quadrature points [33]. The choice of quadrature rule thus has a significant impact on the overall cost of discretization, and count in favor of using few quadrature points when this is feasible. Independent of which strategy is chosen, the local problems can be solved in parallel.

## 3. Multi-Point Methods for thermo-poroelasticity

In this section, we will apply concretely the discretization concepts presented in Section 2 to the problem of thermo-poroelasticity. As in Section 2, our aim is to be pedagogical, and we will defer the discussion of mathematical properties to Section 4.

We will address the discretization for thermo-poroelasticity through four steps, following the natural progression from flow in porous media in Section 3.1, via elasticity in Section 3.2, and then combining the concepts to poroelasticity in Section 3.3. Finally, the full thermo-poroelastic discretization in presented in Section 3.4.

### 3.1 Flow in porous media

The basic equations for flow in porous media, as far as this exposition is concerned, are captured by the (steady state) conservation law for fluid flow, equation (2.1), and Darcy's law relating pressure gradients to fluid flux. In preparation for poroelasticity later, we will denote the pressure potential as $p$, and the fluid flux as $\tau_p$, and re-state conservation and Darcy's law in terms of these variables as

$$\int_{\partial \omega} n \cdot \tau_p \, dS = \int_\omega r_p \, dV \tag{3.1}$$

$$\tau_p = -\kappa \, \nabla p + g \tag{3.2}$$

In a slight abuse of language, will refer to the 2$^{nd}$ order tensor $\kappa$ as the permeability, and $g$ as gravity.

The Multi-Point Flux Approximation (MPFA) follows exactly the general structure of MPxA methods, with equation (3.2) imposed exactly in order to define the discrete constitutive law $\boldsymbol{B}$. We will avoid restating the equations of Section 2.2, and summarize that a numerical fluid normal flux $\boldsymbol{q}$ is defined by the MPFA method as a linear function of pressure, i.e. the finite volume scheme for equations (3.1-3.2) is given as

$$\boldsymbol{D_p q} = \boldsymbol{r_p} \tag{3.3}$$

$$\boldsymbol{q} = \boldsymbol{Q_p p} + \boldsymbol{Q_g g} \tag{3.4}$$

Here we use the subscript $p$ on the discrete divergence operator for this scalar problem, in order to distinguish it from the divergence operator for vector problems in the next sections.

With the choice of the simplified penalty function $\mathcal{M}_s^\eta$ and thus replacing the minimization problem by equation (2.23), the method is referred to simply as the MPFA-O($\eta$), and represents one of the two original MPFA methods [3]. When the full penalty function $\mathcal{M}_s$ is used, we refer to this as the generalized MPFA-O method.

It is not a priori obvious whether we should consider the pressure $\boldsymbol{p}$ as representing a cell-center variable or a mean value for the cell. However, when reviewing the conservation law, we note that the conserved quantity is actually the integrated mass density over the cell, which is related to pressure via a constitutive law. As such, in most implementations, it is most natural to consider the pressure as a mean value for the cell.

For the scalar problem, several specialized variants of the MPFA methods can be derived [3, 21, 22, 23]. However, each of these variants utilize a separate calculation for each of the subface fluxes. As such, these variants cannot be interpreted as having a unique piecewise linear pressure field, and are thus more complex to describe, implement and analyze. Their usage is limited in practice.

## 3.2 Elasticity

The equations of elasticity have the same basic elliptic structure as those for flow, however they have significant differences in the details. First, we note that the deformation vector $u$ takes the role of potential, while the stress tensor $\pi$ takes the role of a flux. The steady state of conservation of momentum is the balance equation for forces

$$\int_{\partial \omega} n \cdot \pi \, dS = \int_\omega r_u \, dV \tag{3.5}$$

Note that this is a vector equation. Elastic materials satisfy Hooke's law, which can be written as

$$\pi = \mathbb{C} : \varepsilon(\nabla u) + \chi \tag{3.6}$$

Here $\varepsilon(\nabla u)$ denotes the material strain, which in the regime of small deformations can be linearized as

$$\varepsilon(\nabla u) = \frac{\nabla u + \nabla u^T}{2} \tag{3.7}$$

The external tensor field $\chi$ can arise from an existing stress state in the material, or as we will see below, from interactions with a separate process in composite materials. We will assume that the external tensor field is always symmetric.

The strain tensor $\varepsilon(\nabla u)$ is symmetric by definition, and we therefore refer to equations (3.6)-(3.7) as Hooke's law with *strong symmetry*. In general, the gradient of the deformation $\nabla u$ need not be symmetric, and as a consequence, the compound action of $\mathbb{C}$ and $\varepsilon$ does not have a unique inverse when Hooke's law is written on the form (3.6)–(3.7). This has consequences for the stability of numerical methods, as we will see below.

We therefore consider also an alternative formulation of Hooke's law, known as Hooke's law with *weak symmetry*. Equations (3.6) and (3.7) can then be equivalently stated as

$$\pi = \mathbb{C} : (\nabla u + b) + \chi \tag{3.8}$$

where $b$ is the asymmetry of the gradient of deformation, which we will at the moment treat as unknown. In order to determine $b$, we enforce that the stress is symmetric. It turns out that it is sufficient to impose symmetry of the stress tensor weakly [34, 33]. Pre-empting that we will impose symmetry on the dual grid, we state the weak symmetry as follows: For all subdomains $\omega^* \in \Omega$, it holds that

$$\int_{\omega^*} as(\pi)\, dS = 0 \tag{3.9}$$

where the asymmetry of a tensor is defined as

$$as(\pi) = \frac{\pi + \pi^T}{2} \tag{3.10}$$

It is straight-forward to verify, by setting $b = -as(\nabla u)$, that equations (3.8)-(3.10) are formally equivalent to equations (3.6)-(3.7), and have also been recently considered in the mixed finite element context (see e.g. [34]).

### 3.2.1 MPSA with strong symmetry

The MPxA finite volume method can be applied directly to equation (3.5-3.7), in exact analogy to the scalar case, and we refer to this as the generalized MPSA-O method. As with the fluid flow, we can use the constitutive law (3.6) directly to define the discrete constitutive law $\boldsymbol{B}$.

Again we will avoid restating the equations of Section 2.2, and summarize that a numerical normal stress (i.e. traction) $\boldsymbol{w}$ is defined by the MPSA method as a linear function of pressure, i.e. the finite volume scheme for equations (3.5-3.7) is given as

$$\boldsymbol{D}_u \boldsymbol{w} = \boldsymbol{r}_u \tag{3.11}$$

$$\boldsymbol{w} = \boldsymbol{W}_u \boldsymbol{u} + \boldsymbol{W}_\chi \boldsymbol{\chi} \tag{3.12}$$

Note that while the action of $\boldsymbol{D}_p$ and $\boldsymbol{D}_u$ are logically similar, the matrixes have slightly different structure as equations (3.11) represent $n$ times as many degrees of freedom due to the vector nature of the elasticity equations.

It turns out that the simplified penalty function $\mathcal{M}_s^\eta$ is not suitable for elasticity with strong symmetry. Indeed, the symmetry of the stress tensor reduces the number of constraints imposed by the local balance stated in equation (2.13), and the minimization problem given by (2.26) for the simplified penalty functions fail to have a unique solution. On the other hand, it can be shown that equation (2.27) does have a solution, and as such, only the generalized MPSA-O method is applicable elasticity with strong symmetry [35]. While the generalized MPSA-O method is well suited for polyhedral grids in 2D and 3D, it is deficient on simplicial meshes [28, 33]. As is the case with mixed finite elements [34], it turns out that the formulation with weakly imposed symmetry is preferable.

### 3.2.2 MPSA with weak symmetry

When we consider the constitutive law with weak symmetry, we quickly note that the MPxA framework needs an adaptation in order to accommodate the condition equation (3.10). Indeed, by imposing equation (3.10) on each dual cell, it is equivalent to stating that for all $s \in \mathcal{V}$, it holds that

$$\sum_{\widetilde{\omega} \in \widetilde{\mathcal{T}}_s} \int_{\widetilde{\omega}} as(\pi) \, dS = 0 \tag{3.13}$$

Since the stress $\pi$ is approximated as constant on each subcell, equation (3.13) can easily be represented in terms of degrees of freedom as the matrix equation

$$\boldsymbol{S\pi} = \boldsymbol{0} \tag{3.14}$$

where again the matrix $\boldsymbol{S}$ can be written as a sum of local matrices for each vertex, e.g. $\boldsymbol{S} = \sum_{s \in \mathcal{V}} \boldsymbol{S}_s$. With this tool in hand, the MPxA framework can be used directly to obtain a discretization, which we refer to as MPSA-W (W signifying weak symmetry), and state to be precise as:

<u>Definition 3.1 (MPSA-W (elasticity)):</u> Let $\mathcal{N}_u$ be the null-space of the vector extension of the constraints given in equation (2.11) – (2.13), as well as (3.14), subject to a given displacement $\boldsymbol{u}$. Then the MPSA-W method is defined by the pair $(\widetilde{\boldsymbol{u}}, \boldsymbol{\pi}) \in \mathcal{N}_u$ such that

$$(\widetilde{\boldsymbol{u}}, \boldsymbol{\pi}) = \arg \min_{(\widetilde{\boldsymbol{u}}', \boldsymbol{\tau}') \in \mathcal{N}_p} \mathcal{M}(\widetilde{\boldsymbol{u}}') \tag{3.16}$$

and the numerical normal stress is defined as $\boldsymbol{w} = \boldsymbol{W}_u \boldsymbol{u} + \boldsymbol{W}_\chi \boldsymbol{\chi} \equiv \boldsymbol{\Sigma}_\mathcal{F} \boldsymbol{F\pi}$.

Clearly, due to the choice of imposing the asymmetry of the stress on the dual grid, the MPSA-W method reduces to local calculations in the same way as other MPxA methods.

The MPSA-W discretization is now obtained by combining the normal stress from the MPSA-W method with the momentum balance, equation (3.11), in exactly the same manner as for the generalized MPSA-O method.

The MPSA-W method behaves qualitatively analogously to the MPFA methods, and can be used together with either the full penalty functions or the simplified penalty functions. In contrast to the generalized MPSA-O method, the MPSA-W method is equally applicable to polyhedral as well as simplicial grids [33].

## 3.3 Poroelasticity

We will consider the linearized equations for poro-elasticity, after an implicit discretization over a time-step length $\theta$. Then the linear system for pressure and displacement consists of two conservation laws for the fluid and solid [27]:

$$\int_\omega \alpha : \nabla u + cp \, dV + \theta \int_{\partial \omega} n \cdot \tau_p \, dS = \int_\omega r_p \, dV \tag{3.17}$$

$$\int_{\partial \omega} n \cdot \pi \, dS = \int_\omega r_u \, dV \tag{3.18}$$

as well as the constitutive laws for fluid flow (Darcy) and stress in poroelastic materials (Biot), stated in the form with weak symmetry:

$$\tau_p = -\kappa \, \nabla p + g \tag{3.19}$$

$$\pi = \mathbb{C} : (\nabla u + b) - \alpha p \tag{3.20}$$

$$\int_{\omega^*} as(\pi) \, dS = 0 \tag{3.21}$$

Relative to the previous sections, we have introduced the Biot coupling coefficient $\alpha$, which is in general a symmetric second-order tensor (but often approximated as a scalar times an isotropic tensor in practice), as well as the effective compressibility term $c$, containing contributions from both bulk and fluid compressibility. Note also that the information from the previous time-step is integrated into the right-hand side term $r_p$.

In order to obtain a numerical stress function for poroelasticity, we follow the MPxA framework outlined in Section 2.2. From the perspective of the mechanics, the pressure is an external stress in the constitutive law, while conversely, from the perspective of flow, the mechanics affects the conservation statement.

We will therefore for the mechanical calculation consider the fluid pressure as the (previously external) imposed stress for the cell. The MPxA framework can then be applied directly, as in the case of Section 3.2. To incorporate the new material constants arising from the coupling terms, let $C$ be the application of the compressibility factor $c$ at the cell level. Moreover, let $A_1$ and $A_2$ be the application of the Biot coefficient $\alpha$ at the subcell level, where $A_1$ acts as the double inner-product on tensors (confer equation (3.17)), while $A_2$ acts a tensor-scalar product (confer equation (3.20)). Finally, similarly to $\Sigma_\mathcal{F}$ we let $\Sigma_\mathcal{T}$ be the summation over subcells, weighted by the subcell volumes.

With these definitions, the MPSA method can be directly applied to the linearized equations for poroelasticity. For completeness, we state its definition as:

<u>Definition 3.2 (MPSA-W (poroelasticity))</u>: Let $\mathcal{N}_{u,p}$ be the null-space of the vector extension of the constraints given in equation (2.11) – (2.13), as well as (3.14), subject to a given displacement $u$ and a pressure $p$. Then the MPSA-W method for poroelasticity is defined by the pair $(\tilde{u}, \pi) \in \mathcal{N}_{u,p}$ such that

$$(\tilde{u}, \pi) = \arg \min_{(\tilde{u}', \pi') \in \mathcal{N}_{u,p}} \mathcal{M}(\tilde{u}') \tag{3.24}$$

and the numerical normal stress is defined as $w = W(u, p) = W_u u + W_p p \equiv \Sigma_\mathcal{F} F \pi - \Sigma_\mathcal{F} F A_2 E^* p$. Moreover, the impact of displacement on the fluid mass conservation law is given by $J(u, p) = J_u u + J_p p \equiv \Sigma_\mathcal{F} A_1 G \tilde{u}$.

We note that the linear discretization matrices $W_u, W_p, J_u$ and $J_p$ are implicitly defined, since $\tilde{u}$ and $\pi$ are linear functions of $u$ and $p$. As previously, all the minimization problems can be solved in parallel for each vertex of the grid. Moreover, the same points about simplified penalty functions as discussed in Section 3.2.2 are applicable.

We close this section by stating the full discrete system for the poroelastic equations (3.17-3.21):

$$J_u u + J_p p + Cp + \theta D_p q = r_p \tag{3.25}$$

$$D_u w = r_u \tag{3.26}$$

$$q = Q_p p + Q_g g \tag{3.28}$$

$$w = W_u u + W_p p \tag{3.29}$$

By eliminating the flux and normal stress, we get a system of matrix equations only in terms of cell-center pressure and displacement, given as

$$\begin{pmatrix} D_u W_u & D W_p \\ J_u u & C + \theta D_p Q_p + J_p \end{pmatrix} \begin{pmatrix} u \\ p \end{pmatrix} = \begin{pmatrix} r_u \\ r_p - \theta D_p Q_g g \end{pmatrix} \tag{3.30}$$

As expected, equation (3.25-3.29), and consequently also (3.30), has essentially the same structure as the continuous problem. The exception is the presence of the term $J_p$, which appears implicitly in the pressure conservation equation, since the discrete displacement gradient $\nabla u$ calculated by the MPSA method also depends on the pressure. This dependence is weak and can be interpreted as the expansion of a (sub)cell due to an increment of pressure in that cell. As such, it has the structure of a Laplacian operator, scaled by the bulk modulus of the solid $\mathbb{C}$, and the square of the characteristic length scale of the cell $\Delta x$, i.e., the following spectral equivalence holds [27]:

$$J_p \sim \frac{(\Delta x)^2}{\mathbb{C}} D_p D_p^T \tag{3.31}$$

## 3.4 Thermo-Poroelasticity

As the final application of the MPxA framework, we will consider the transport of heat in a poroelastic medium. Again, we will consider the equations subject to an implicit discretization over a time-step length $\theta$. One has a choice in which variable to use to represent heat, however we will for simplicity consider temperature $\phi$, as we are primarily interested in the spatial discretization of the linearized equations. Then the linear system for pressure, displacement and temperature then consists of three conservation laws for the fluid, entropy, and solid [8, 36]:

$$\int_\omega \alpha_p : \nabla u + c_{p,p} p + c_{p,\phi} \phi \ dV + \theta \int_{\partial \omega} n \cdot \tau_p \ dS = \int_\omega r_p \ dV \tag{3.32}$$

$$\int_\omega \alpha_\phi : \nabla u + c_{\phi,p} p + c_{\phi,\phi} \phi \ dV + \theta \int_{\partial \omega} n \cdot \tau_\phi \ dS = \int_\omega r_\phi \ dV \tag{3.33}$$

$$\int_{\partial \omega} n \cdot \pi \ dS = \int_\omega r_u \ dV \tag{3.34}$$

as well as the constitutive laws for fluid flow (Darcy), heat transfer (Fourier's law and advection), and stress in thermo-poroelastic materials, the latter stated in the form with weak symmetry:

$$\tau_p = -\kappa_p \nabla p + g \tag{3.35}$$

$$\tau_\phi = -\kappa_\phi \nabla \phi + \phi \tau_p \tag{3.36}$$

$$\pi = \mathbb{C} : (\nabla u + b) - \alpha_p p - \alpha_\phi \phi \tag{3.37}$$

$$\int_{\omega^*} as(\pi) \, dS = 0 \tag{3.38}$$

Relative to the previous sections, we have for each of the fluid and thermal conservation laws separate linearized constitutive laws $c_{\phi,p}$, $c_{\phi,\phi}$, $c_{p,p}$, and $c_{p,\phi}$, Biot coupling coefficients $\alpha_p$ and $\alpha_\phi$, and constitutive laws $\kappa_p$ and $\kappa_\phi$. Furthermore, we notice that the constitutive law for heat flux, given in equation (3.36), is non-linear due to the presence of the product $\phi \tau_p$, representing heat advection with the fluid flux.

With exception of the head advection term, the coupled problem represented by equations (3.32-3.38) presents no new challenges relative to the poroelastic problem considered in section 3.3, and the application of the MPxA method for the problem is equivalent. Thus we have the following discrete system for thermo-poroelasticity, which is the discrete analog of equations (3.32-3.38). The conservation laws take the form:

$$J_{p,u} u + J_{p,p} p + C_{p,p} p + C_{p,\phi} \phi + \theta D_p q_p = r_p \tag{3.39}$$

$$J_{\phi,u} u + J_{\phi,\phi} \phi + C_{\phi,p} p + C_{\phi,\phi} \phi + \theta D_p q_\phi = r_\phi \tag{3.40}$$

$$D_u w = r_u \tag{3.41}$$

While the constitutive laws take the form:

$$q_p = Q_{p,p} p + Q_{p,g} g \tag{3.42}$$

$$q_\phi = Q_{\phi,\phi} \phi + \phi^* q_p \tag{3.43}$$

$$w = W_u u + W_p p + W_\phi \phi \tag{3.44}$$

The discrete matrixes are constructed exactly as in Definition 3.2, with the notational convention that (say) $W_\phi$ is calculated using the coupling coefficient $\alpha_\phi$, while $W_p$ is calculated using the coupling coefficient $\alpha_p$. Similarly, the discrete flux stencil $Q_{\phi,\phi}$ is calculated using the MPFA method of Section 3.1, with the coefficient tensor $\kappa_\phi$.

It remains to define the temperature $\phi^*$ on faces of the grid. We denote the matrix with these temperatures on the main diagonal as $\phi^*$, for which the simplest and most commonly choice is obtained via the so-called upstream weighting [37, 7]. Thus, for any face $\sigma \in \mathcal{F}$, with neighboring cells $\omega_{k_1}$ and $\omega_{k_2}$ where $k_1 < k_2$, then

$$\phi^*_{\sigma,\sigma} = \begin{cases} \phi_{k_1} & \text{if } q_{p,\sigma} \geq 0 \\ \phi_{k_2} & \text{if } q_{p,\sigma} < 0 \end{cases} \tag{3.45}$$

The entries of $\phi^*$ are typically taken as zero away from the main diagonal.

A compact and simple discretization for coupled thermo-poromechanics is then obtained in terms of cell-center variables as

$$\begin{pmatrix} D_u W_u & D_u W_p & D_u W_\phi \\ J_{p,u} u & C_{p,p} + \theta D_p Q_{p,p} + J_{p,p} & C_{p,\phi} \\ J_{\phi,u} u & C_{\phi,p} + \theta D_p \phi^* Q_{p,p} & C_{\phi,\phi} + \theta D_p Q_{\phi,\phi} + J_{\phi,\phi} \end{pmatrix} \begin{pmatrix} u \\ p \\ \phi \end{pmatrix} = \begin{pmatrix} r_u \\ r_p - \theta D_p Q_{p,g} g \\ r_\phi - \theta D_p \phi^* Q_{p,g} g \end{pmatrix} \quad (3.30)$$

Note that this discretization inherits the non-linearity of the original problem (3.32-3.38), due to the presence of the advective term, which explicitly becomes $\theta D_p \phi^* Q_{p,p} p$.

# 4. Mathematical properties of MPxA methods

Since their inception, the MPFA (and later MPSA) methods have been intensely studied. Various viewpoints have been considered, using both analysis frameworks building on theory of finite volume and mixed finite element methods, as well as numerical validations. As a whole, these studies provide a comprehensive perspective on not just the properties of the MPFA finite volume discretization for the model problem from the introduction, equations (1.2-1.2), but also the performance for elasticity and coupled problems as discussed in Sections 3.2-3.4. We will summarize some of the main aspects below.

## 4.1 Analysis of consistency and convergence

Already in the earliest papers on MPFA methods, the consistency of the discretization was validated on parallelogram grids [3, 4]. More general analysis followed a decade later, and the first proof of convergence was established by Klausen and Winther, considering perturbations of parallelogram grids [38]. That analysis explicitly constructed a mixed finite element method which is algebraically equivalent to the MPFA method with simplified quadrature, by using the local problems detailed in Section 2.2 and 2.3 to define so-called "broken" finite element spaces for the flux.

While the analysis of Klausen, Winther provides both convergence as well as rates of convergence, it is quite restrictive, and does not apply to polyhedral grids nor non-smooth coefficients. The analysis was later extended to general polyhedral grids by Klausen and Stephansen by exploiting a link to mimetic finite differences [39]. A different approach was therefore pursued by Agelas et al [26], where they considered a formulation of the MPFA method in terms of the finite volume analysis framework [24]. This yielded convergence proofs through compactness arguments for quite general grids, and with minimal assumptions on the coefficients. On the other hand, this generality reduces the regularity of the exact solution, and thus explicit rates of convergence cannot be considered in this framework. Furthermore, the proofs suffered from an a priori assumption that the local formulation of MPFA was uniformly coercive (with respect to all corners of the grid and all grid refinement). Such an assumption automatically holds for self-similar grid refinements, but was not proved.

The approach of Agelas was extended to show the convergence of MPSA for elasticity, and later also to show the convergence of MPFA+MPSA for the poroelastic problem of section 3.3 [28, 27]. In these proofs, the general case of penalty functions with multiple quadrature points, as introduced in Section 2.2.3, was first considered. Considering the full penalty formulation had the further advantage of avoiding the local coercivity assumptions of Agelas, as the local coercivity could be proved based on the structure of the minimization problems. Moreover, the convergence proofs were shown to hold even for degenerate coefficients, such as incompressible materials and near-zero time-step size.

It is worth noting the related development of so-called Multipoint Flux Mixed finite Element (MFME) [40] and Multipoint Stress Mixed Finite Element (MFSE) [41] methods. These methods are obtained from mixed-finite element methods with BDM1 elements for flux (or stress) and P0 elements for pressure (or displacement), using various quadrature rules to eliminate the flux variables. The resulting methods, for which detailed analysis is possible based on standard theory of mixed finite elements, are close cousins of the MPxA finite volume methods described herein. However, these methods are less suited for geometrically complex problems, since the quadrature rules lead to reduced rates of convergence for rough grids [42], and the underlying finite element spaces preclude the applications to polyhedral grids.

### 4.2 Monotonicity

The question of monotonicity is essentially a translation of Hopf's lemma from the continuous problem to the discrete problem. For the scalar case, Hopf's lemma can be stated as the property that for a zero right-hand side $r_p = 0$, then the maximum (and minimum) value of the solution $p$ should be found on the boundary of the domain [43].

Several numerical methods preserve the monotonicity property, in particular those that lead to discretization matrices on the form of $M$-matrices such as TPFA and FE. On the other hand, this property is in no way guaranteed, and as an example, the MFE method is in general not monotone. As the MPFA discretization does not guarantee an $M$-matrix, the question of monotonicity is subtle.

Sufficient and necessary conditions for any finite volume discretization to satisfy a discrete maximum principle can be established in the case of quadrilateral grids [29, 44]. As a result, it is now known that there are essentially three categories of grid (cells): 1) Those for which essentially any finite volume methods will lead to monotone discretization, 2) Those for which it is possible by a judicious choice to construct a monotone discretization, and 3) Those for which no linear finite volume discretization (with a relatively compact stencil) exists.

Clearly, point 3) above means that there are certain grids which are sufficiently bad that the performance of a MPxA discretization cannot be guaranteed, and these are in general grids combining high aspect ratios with a high degree of skewness. Point 2) above furthermore inspired research into constructing MPxA methods that are optimal with respect to monotonicity. Such methods can be constructed either by optimizing the location of quadrature points in the penalty functions [45, 29], or by allowing for more general formulation of the MPxA methods than that outlined in Section 2.2. As a result of the latter approach, the MPxA-Z method with a larger stencil [22], and the MPxA-L method with a smaller stencil [23, 21], were developed.

### 4.3 Numerical investigations of convergence

In addition to the analysis summarized above, it is worth noting that the convergence properties of the MPFA and MPSA methods have been extensively studied numerically. These numerical investigations also consider problems not covered by analysis, due to either challenging coefficients [46], non-linearities [30], or grids [47]. We will review some of these results here, emphasizing the results that give a most comprehensive understanding of the general features of the MPxA methods.

### 4.3.1 Convergence rates for smooth solutions

For problems with smooth coefficients on regular domains, the analysis of MPFA methods indicates that one can expect 2nd order convergence of the potential and 1st order convergence of the fluxes [38]. In practice 2nd order convergence of fluxes has been observed in numerical calculations for the flow problem, and what appears to be 1.5 order convergence for the elasticity and Biot problems. We will revisit some of these results here [27].

The problem under consideration is the poroelastic equations as presented in Section 3.3, with the MPFA and MPSA methods using full penalty functions as given in equation (2.20), and with the elasticity discretized with strong symmetry.

With the $L^2$ norms defined as

$$\|u\|_{\mathcal{T},0} = \left(\sum_{k\in\mathcal{T}} m_k u_k^2\right)^{1/2} \quad \text{and} \quad \|q\|_{\mathcal{F},0} = \left(\sum_{\sigma\in\mathcal{F}} m_\sigma^2 q_\sigma^2\right)^{1/2} \quad (4.1)$$

We can define errors using the following $L^2$ type metrics, where variables in plain type are the exact analytical solution, and variables in bold are the discrete solutions, as in the preceding sections. The error in primary variables is then measured as relative to the projection $\Pi_\mathcal{T}$ which returns cell-center values (i.e. $(\Pi_\mathcal{T} p)_k = p(x_k)$):

$$\epsilon_u = \frac{\|\boldsymbol{u}-\Pi_\mathcal{T} u\|_{\mathcal{T},0}}{\|\Pi_\mathcal{T} u\|_{\mathcal{T},0}} \quad \text{and} \quad \epsilon_p = \frac{\|\boldsymbol{p}-\Pi_\mathcal{T} p\|_{\mathcal{T},0}}{\|\Pi_\mathcal{T} p\|_{\mathcal{T},0}} \quad (4.2)$$

and the error in secondary variables as relative to the projection $\Pi_\mathcal{F}$ which returns face-center fluxes (i.e. $(\Pi_\mathcal{F}\tau)_\sigma = \tau(x_\sigma)\cdot n_\sigma$):

$$\epsilon_\pi = \frac{\|\boldsymbol{w}-\Pi_\mathcal{F} u\|_{\mathcal{T},0}}{\|\Pi_\mathcal{F} u\|_{\mathcal{T},0}} \quad \text{and} \quad \epsilon_q = \frac{\|\boldsymbol{q}-\Pi_\mathcal{F}\tau_p\|_{\mathcal{T},0}}{\|\Pi_\mathcal{F}\tau_p\|_{\mathcal{T},0}} \quad (4.3)$$

Finally, we also consider the error based on the $L^2$ seminorm of pressure, which discards the datum value, defined as

$$\epsilon_{p,|} = \inf_{p_0\in\mathbb{R}} \frac{\|\boldsymbol{p}-\Pi_\mathcal{T} p + p_0\|_{\mathcal{T},0}}{\|\Pi_\mathcal{T} p\|_{\mathcal{T},0}} \quad (4.4)$$

In order to illustrate the numerical convergence rate of the primary variables, we give the *primary error* associated with the primary variables displacement and pressure as

$$\epsilon_{\boldsymbol{u},p} = \epsilon_{\boldsymbol{u}} + c\epsilon_p \quad (4.5)$$

Furthermore, it can be shown that the MPSA-MPFA discretization is stable even for degenerate timestep size $\tau \to 0$ and compressibility $c \to 0$, subject to a weighted combination of the norms above [27]. Thus, we introduce the so-called *stable error*

$$\epsilon_\Sigma = \epsilon_{\boldsymbol{u}} + \epsilon_\pi + (\theta + c)\epsilon_p + \theta\epsilon_q + \epsilon_{p,|} \quad (4.6)$$

The numerical convergence rates for a smooth manufactured solution on irregular simplicial, quadrilateral, and polyhedral grids, as illustrated in Figure 4.1. The calculations are based on seven levels of refinement for each grid type, for which the finest grid level has a characteristic cell diameter of $h \sim 2^{-7}$, the results of which are summarized in Table 4.1. We note that as expected, 2nd order convergence is observed for primary variables, and better-than-1st order convergence is observed for fluxes and stresses.

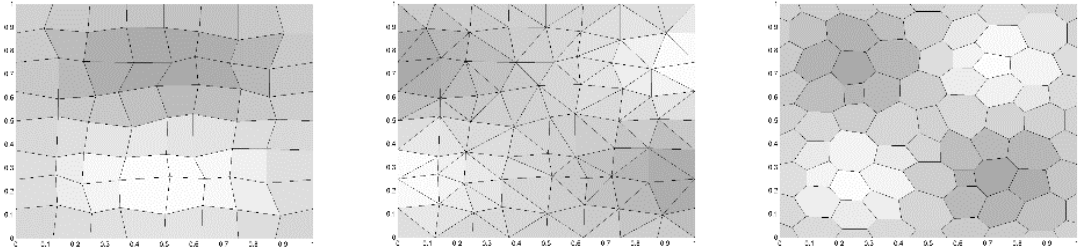

**Figures 1**: From left to right the figures illustrate grid types A (quadrilaterals), B (triangles), C (unstructured grids). Furthermore, in grey-scale, the figures indicate the structure of the analytical solution used for this example.

| $\epsilon_{u,p}$ | $\tau = 1$ | | | $\tau = 10^{-6}$ | | |
|---|---|---|---|---|---|---|
| Grid | A | B | C | A | B | C |
| $c = 1$ | 2.00 | 1.97 | 1.99 | 1.99 | 1.95 | 1.98 |
| $c = 10^{-2}$ | 2.00 | 1.97 | 1.99 | 1.98 | 1.94 | 1.98 |
| $c = 10^{-6}$ | 2.00 | 1.97 | 1.99 | 1.98 | 1.94 | 1.98 |

| $\epsilon_\Sigma$ | $\tau = 1$ | | | $\tau = 10^{-6}$ | | |
|---|---|---|---|---|---|---|
| Grid | A | B | C | A | B | C |
| $c = 1$ | 1.36 | 1.36 | 1.27 | 1.09 | 1.14 | 1.16 |
| $c = 10^{-2}$ | 1.32 | 1.32 | 1.23 | 1.20 | 1.29 | 1.28 |
| $c = 10^{-6}$ | 1.32 | 1.32 | 1.23 | 1.20 | 1.29 | 1.29 |

**Tables 1**: Asymptotic convergence rate of stable primary error $\epsilon_{u,p}$ and stable error $\epsilon_\Sigma$ for grids of types A, B, and C.

### 4.3.2 Convergence rates for singular solutions

In order to assess the convergence rates for non-smooth problems, Eigestad and Klausen considered domains with discontinuous permeability coefficients, such as illustrated in figure 4.2 [46].

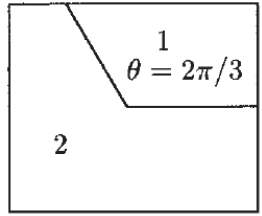

**Figure 4.2**: Partitioning of domain such that a non-trivial material discontinuity can be defined.

For such domains, analytical solutions can be defined on using polar coordinates around the center point, on the form

$$p(r,\theta) = r^\alpha (a_i \cos(\alpha\theta) + b_i \sin(\alpha\theta)) \tag{4.7}$$

The constants $\alpha$, $a_i$ and $b_i$, for $i = 1,2$, depend on the permeability contrast chosen, and in particular, the exponent $\alpha$ also determines the regularity of the solution, i.e.

$$p \in H^{1+\alpha} \tag{4.8}$$

For such problems, they report a loss of convergence rate, such that one observes that the pressure converges at a rate of $\epsilon_p \sim h^{\min(2,2\alpha)}$ while the flux converges at a rate $\epsilon_q \sim h^{\min(1,\alpha)}$. For the particular choice of $\theta = 2\pi/3$, and a permeability contrast $\frac{k_1}{k_2} = 100$, the resulting analytical solution has the exponent $\alpha \approx 0.75$. Figure 4.3 illustrates the convergence for this case, based on the MPFA method with simplified penalty functions ($\eta = 0$), and both regular and perturbed grid sequences.

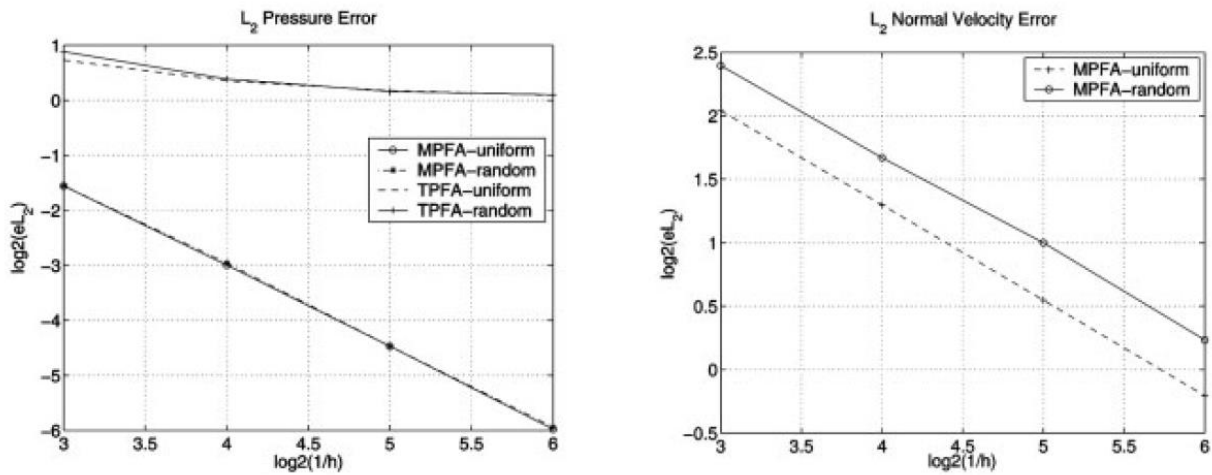

**Figure 4.3**: Convergence of pressure (left) and flux (right) for the MPFA method for the non-smooth problem of Section 4.3.1. The observed convergence rates for this problem, where the exact solution $p \in H^{1.5}$ is order 1.5 for pressure and order 0.75 for flux. For comparison, the TPFA method was included in the study, which does not converge. Figure reused, with permission, from [46].

### 4.3.3 Robustness on degenerate grids

Contrasting the previous two studies, Nilsen et al. emphasized degeneracies of the grid (as opposed to regularity-preserving refinements) [47]. To this end, they considered a series of cases with polyhedral grids, grids of high aspect ratio, and unusual refinement strategies. All of their calculations considered the MPSA discretization with strong symmetry, applied to either elasticity or coupled with MPFA for Biot.

An illustrative example from that study, considers a problem of non-matching grids, meeting at a thin layer, as illustrated in Figure 4.4 (left). The thin layer is discretized by a finer grid which has roughly isotropic shape, as shown in Figure 4.4 (right). The color scale in that figure indicates the approximation error relative to a smooth reference solution, which can be seen to be less than 4% in displacement (left figure) and as much as around 50% in the grid cells immediately adjacent to the thin layer for the volumetric strain (right figure).

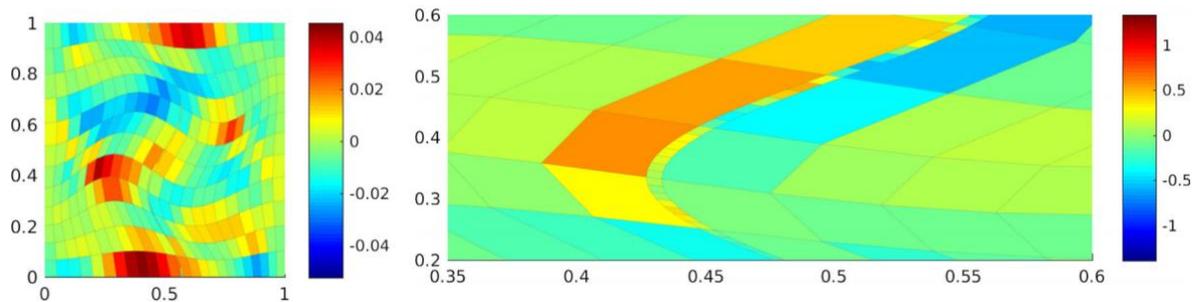

**Figure 4.4**: Illustration of the grid used for robustness study (left). Note there is a vertical section of thin cells in the middle of the domain, mimicking a thin geological layer, as shown in the zoom (right). The discretization on the right-hand side of the thin layer is intentionally chosen to be slightly coarser than the left-hand side to ensure that the inner layer of grid cells is always non-matching relative to the surroundings. The color map on the left indicates the relative error in the x-component of displacement, while the color map on the right indicates the relative error in volumetric strain, both as compared to a manufactured analytical solution for this problem.

To study the robustness of the method the ratio of the thin layer as compared to the external grid cells was varied from a factor 1 to a factor 20 (for comparison, figure 4.4. illustrates a factor 7 difference in grids). Recall that the grid cells in the thin layer are nearly isotropic in shape, thus when the thickness of the thin layer is reduced, the number of cells in the layer is simultaneously increased, thus introducing an increasing number of hanging nodes between. The study can thus be seen both as a study of robustness to an abrupt change in grid sizes in the discretization, as well as a study in the robustness to hanging nodes.

The results are shown in Figures 4.5, where a comparison is also made to a Virtual Element Discretization for the same grid [48]. As can be seen, the approximation quality of the MPSA method is essentially unaffected by the presence of the thin layer, both in the $L^2$ and $L^\infty$ norms. In particular, the stability in the $L^\infty$ norm shows that spurious oscillations are not introduced in the transition between the grids.

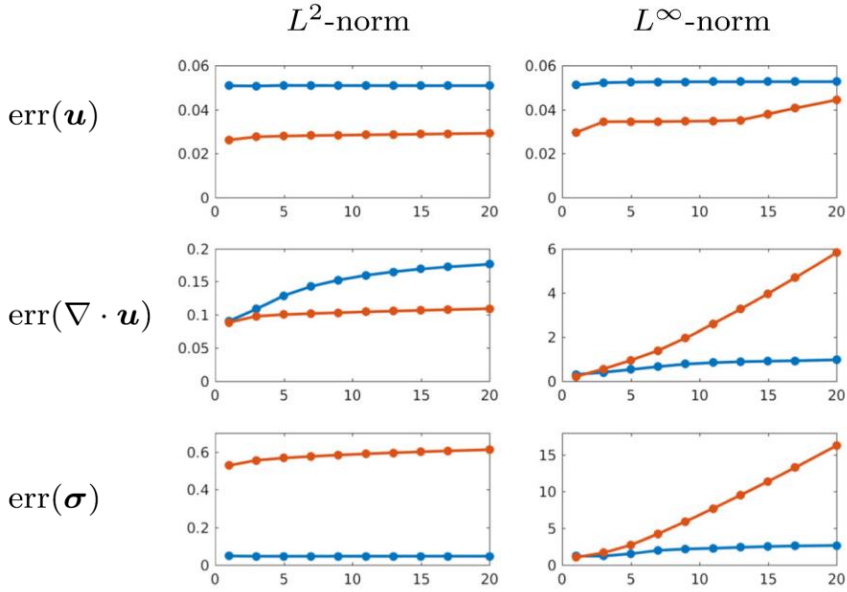

**Figure 4.5**: Error in numerical approximation relative to analytical solution for grid types as shown in Figure 4.5. Blue lines are the MPSA method with strong symmetry from section 3.2.1, while for comparison, an elasticity discretization using the virtual element method [49] is also shown. Errors are shown for displacement, volumetric strain, and stress in the rows, respectively, and using the $L^2$ and maximum norms in the columns. The *x*-axis of all figures denotes the aspect ratio between the thin layer and the outer grid, thus the right-most data-point corresponds to a factor 20 finer grid in the thin layer.

### 4.3.4   Convergence for Thermo-poroelasticity

We close this section with a convergence study for the MPxA discretization of the full thermo-poroelastic problem. To our knowledge, results for this problem have not been reported before. The domain is the unit square, and the grid is formed by quadrilaterals that are roughly perturbed on all refinement levels. As in the study discussed in section 4.3.1, we consider a single time step for the system, with time discretization by a backward Euler approach. All parameters are assigned unit values in this case.  The manufactured solution is given by

$$u = \begin{pmatrix} \sin(2\pi x)\, y(1-y) \\ \sin(2\pi x) \sin(2\pi y) \end{pmatrix}, \quad p = \sin(2\pi x) y(1-y), \quad \phi = xy(1-x)(1-y)$$

The thermo-poroelastic system was discussed with MPSA/MPFA as discussed in Section 3.4, while a single point upstream approach was applied for the temperature advection term.

The convergence behavior is shown in Figure 4.6. Displacement and pressure retain the second order convergence observed on the comparable test for the poro-elastic system considered in Section 4.3.1. For the temperature, the first order scheme for advection makes the convergence deteriorate to first order as the grid is refined, as expected from the theory of hyperbolic conservation laws [37]. Without the advective term, temperature also showed second order convergence. Finally, the mechanical stress and the fluid fluxes both are first order convergent on the perturbed grids. The test case thus confirms that the combined MPxA schemes can be applied successfully applied to problems including a non-linear advection term.

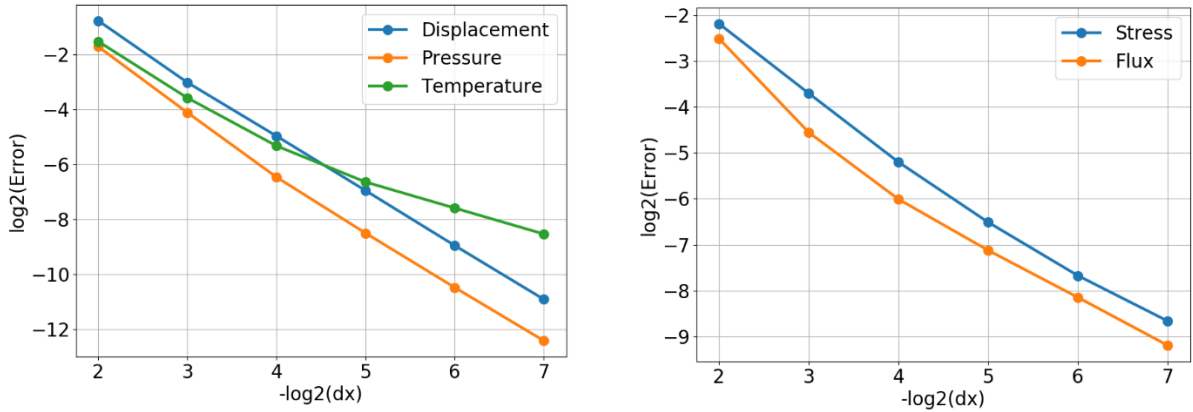

**Figure 4.6:** Convergence plot for the thermo-poroelasticity problem described in section 4.3.4. The convergence for the primary variables is shown to the left, to the right is shown the convergence results for mechanical stresses and fluid fluxes.

## 5. Applications to complex problems

Having reviewed the mathematical properties of MPxA methods in the previous section, we here show application-motivated scenarios where the methods can be applied. We present three setups: Poro-elastic deformation during fluid injection, thermo-poroelastic response to cooling, and flow through a fractured porous media. The cases are designed to showcase the applicability of MPxA methods on a wide range of grids, and the cells is the three cases are respectively perturbed hexahedra, prismatically extended polygons and simplexes. All simulations use the open source simulation tool PorePy [50], which provides an implementation of MPxA method that follow the principles discussed in this chapter, see [51, 50] for details.

### 5.1 Poro-elastic response to fluid injection

We consider the poro-elastic response to fluid injection a domain of $10 \times 10 \times 1.8$ km, covered by $71 \times 71 \times 40 = 201640$ cells forming a Cartesian grid. The test case is motivated by $CO_2$ storage, although only single-phase flow is considered, with alternating layers of high and low-permeable domains that act as storage formation and trap, respectively [5]. The permeability contrast is four to five orders of magnitude, while the elastic moduli are heterogeneous, though of comparable size. The height of the layers varies, so that the computational cells are perturbed from their original hexahedral form, as indicated in Figure 5.1.

The system is discretized with MPSA/MPFA, and fluid injection in the middle of storage layer was simulated. Figure 5.1 shows the fluid pressure and the vertical displacement in a cut domain. The pressure solution adapts to the permeability contrast, while the displacement various smoothly throughout the domain.

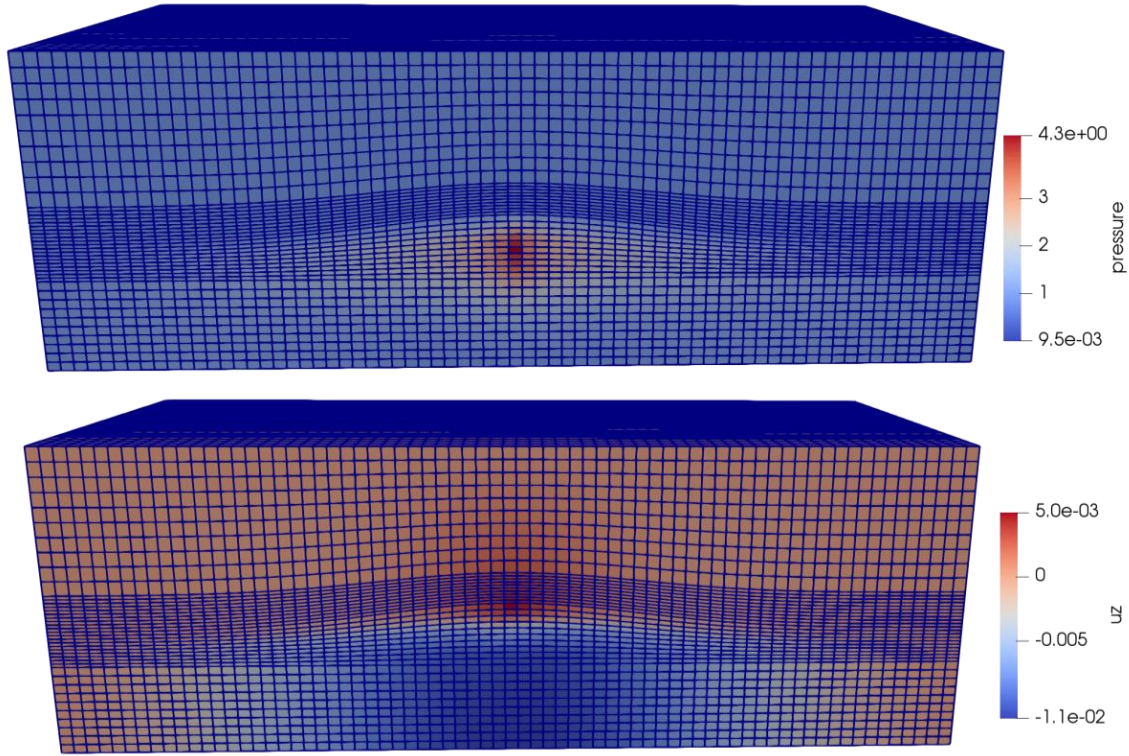

**Figure 5.1**: Simulation of fluid injection into a poro-elastic cube. Fluid is injected into the middle of the domain, and the domain is cut to show the effects near the injection point. Top: Increase in fluid pressure due to injection, measured in bar. Bottom: Vertical displacement, measured in meters.

## 5.2 Thermo-poroelastic response to cooling

To illustrate MPxA applied to the full thermo-poroelastic system, we consider a 3D unit cube that undergoes cooling. Specifically: The domain is cooled at the bottom by a fixing a temperature lower than the initial state. Fluid is allowed to leave through the top, the bottom is impermeable for fluid flow, while the lateral sides are assigned homogeneous Neumann conditions for both fluid and temperature. The domain is fixed on all sides except the bottom, which is free to move. The domain is meshed with polyhedral cells formed by first taking the Voronoi diagram of a 2d triangulation, and then extruding the grid in the third direction. The resulting grid has 4275 cells, with a mixture of 6, 7 and 8 faces per cell.

On this mesh, the full thermo-poroelastic system is discretized as described in Section 3.4. Snapshots of the time evolution of temperature, pressure and displacement are shown in Figure 5.2. At an early stage, the couplings in the system lead to noticeable 3d effects towards the bottom of the domain and significant displacements. The snapshots at later stages reflect the gradual cooling of the domain, and a decrease in pressure and displacement gradients. Note also that the pressure has low regularity at early time, as is expected since the elliptic term for the pressure in equation (3.30) scales with $\theta$. This is consistent with the use of a weighted norm in equation (4.6).

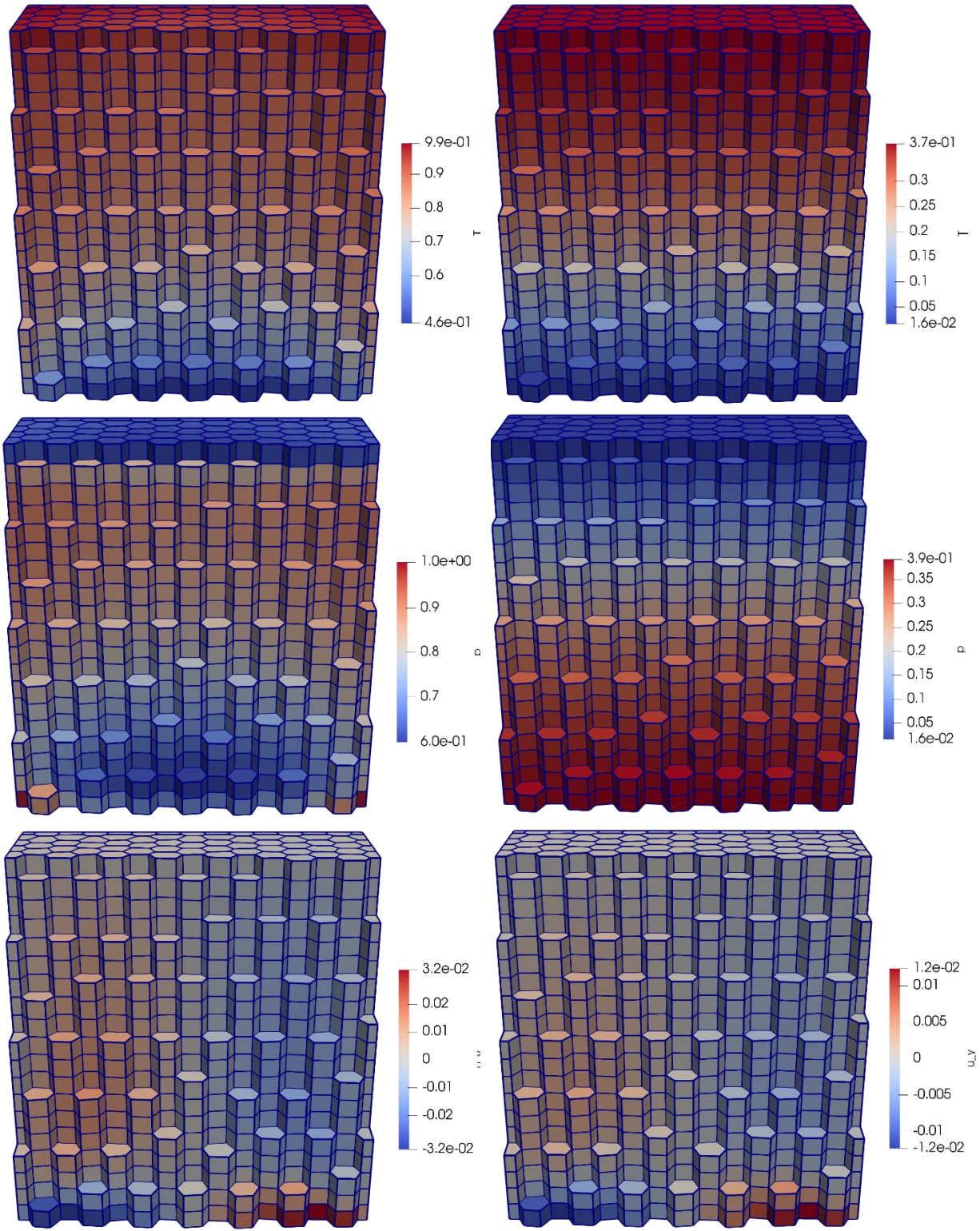

**Figure 5.2:** Thermo-poroelastic deformation of a domain meshed by a polyhedral grid. The domain is cut to expose the 3d structure of the grid cells. The figure shows temperature (top), pressure (middle) and displacement in a direction approximately parallel to the cutting plane (bottom) at an early (left) and late (right) stage.

## 5.3 Flow in fractured porous media

Our final example considers simulation of flow in a 3d domain that contains a network of intersecting fractures. The fractures are modeled as manifolds of co-dimension 1 that are embedded in the host medium. Intersections between fractures form lines of co-dimension 2, while the intersection of intersection lines define intersection points. The fracture network and its host medium thus together define a hierarchy of domains with decreasing dimensions, which we refer to as a mixed-dimensional geometry. Following the model defined e.g. in [51] flow in each of the subdomains is modeled by equations (1.1) – (1.2), with the modification that (1.2) is void in 0d domains.

To define the coupling between subdomains, let $\Omega_h$ and $\Omega_l$ be two domains so that a part of the boundary $\partial \Omega_h$ geometrically coincides with $\Omega_l$, and let $\Gamma$ be an interface between the subdomains. The flow over $\Gamma$ is then governed by the Darcy-like flux law

$$\lambda = \kappa(tr\, p_h - p_l)$$

Here, $\lambda$ is the interface flux, $\kappa$ is the interface permeability, $tr\, p_h$ denotes the trace of the pressure in $\Omega_h$, evaluated on the relevant part of $\partial \Omega_h$, and $p_l$ is the pressure in $\Omega_l$. The interface flux can be considered a mortar variable, which is represented as a Neumann boundary condition to $\Omega_h$ and a source term for $\Omega_l$.

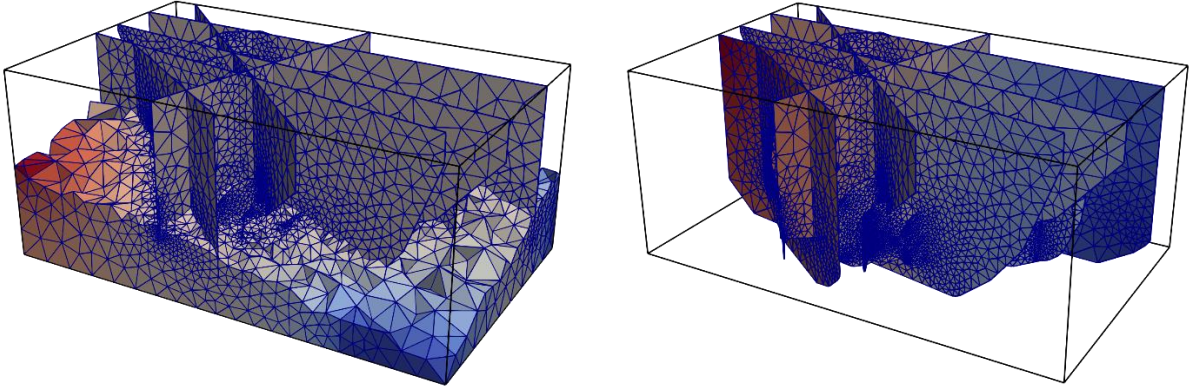

**Figure 5.3:** MPFA pressure solution for the fracture flow benchmark problem presented in Section 5.3. The left figure shows the pressure in the host medium and the fracture network, while the pressure variations internal to the fracture network are illustrated in the right figure.

Following the principles outlined in [51], the MPFA discretization can readily be adapted to mixed-dimensional flow problems. As an illustration we consider the final test case in the benchmark study proposed in [52]. The case consists of a 3d domain with 52 fractures that further form 106 intersection lines as indicated in Figure 5.3. Boundary conditions are set up to drive flow through the host domain and the fracture network, with an inlet in the upper left corner referring to Figure 5.3, and outlets in the two corners of the domain that are in the lower left part of the figure.

The computational mesh is constructed to conform with all fractures and fracture intersection lines. The resulting mesh consists of almost 260K 3d cells, 52k 2d cells (on fracture planes), 1.6k 1d cells (intersection lines) and 105k mortar variables. The MPFA discretization of the full problem produce almost 23M degrees of freedom. The resulting pressure profile is shown in Figure 5.3. The test case thus illustrates the applicability of the MPFA method also to non-standard problems such as flow in mixed-dimensional geometries.

# Acknowledgements

We wish to thank Ivar Aavatsmark, Inga Berre, Geir Terje Eigestad and Ivar Stefansson, for fruitful discussions, and for contributions to various aspects of the numerical implementation. This work forms part of Norwegian Research Council projects 250223 and 267908.